\documentclass[12pt]{article}

\usepackage{setspace}
\usepackage[margin=1in]{geometry}
\usepackage{amsmath}
\usepackage{amssymb}
\usepackage{amsthm}
\usepackage{graphicx,multirow}
\usepackage{algorithm}
\usepackage{algorithmic}
\usepackage{float}
\usepackage{subfigure}
\usepackage{multicol,multirow,float}
\usepackage{bbm}
\usepackage{url}
\usepackage{mathtools}
\usepackage{tabu}
\usepackage{tikz,pgfplots}
\usetikzlibrary{matrix}
\usepgfplotslibrary{groupplots}
\usepackage{authblk}
\usepackage[round,sort]{natbib}
\usepackage[pdfborder={0 0 0},breaklinks=true]{hyperref}
\usepackage{breakcites}
\newtheorem{theorem}{Theorem}

\newtheorem{corollary}[theorem]{Corollary}

\newtheorem{proposition}[theorem]{Proposition}

\newtheorem{remark}{Remark}

\newcommand{\argmin}{\operatornamewithlimits{argmin}}

\newcommand{\ind}{\mathbbm{1}}

\def\ZZ{\mathbb{Z}}
\def\QQ{\mathbb{Q}}

\def\RR{\mathbb{R}}
\def\NN{\mathbb{N}}
\def\EE{\mathbb{E}}
\def\PP{\mathbb{P}}

\def\T{{\cal T}}

\def\B{{\cal B}}

\def\F{{\cal F}}

\def\M{{\cal M}}

\def\F{{\cal F}}

\DeclarePairedDelimiter{\norm}{\lVert}{\rVert}
\def\st{{\rm s.t.}}

\DeclareMathOperator*{\dist}{dist}

\def\st{{\rm s.t.}}
\newcommand{\Limsup}{\operatornamewithlimits{Lim\text{ }sup}}
\newcommand{\co}{\operatorname{co}}
\DeclareMathOperator*{\sgn}{sgn}
\newcommand{\mean}{\operatorname{mean}}

\begin{document}

%\title{Sparse Training with Lipschitz Continuous Loss\\ Functions and a Weighted Group $l_0$-norm Constraint}
%
%\author{\name Michael R. Metel \email michael.metel@huawei.com\\
%       \addr Huawei Noah's Ark Lab\\
%       Montr\'eal, Qc, Canada}
%
%\editor{}
%
%\maketitle	

\onecolumn
\title{Sparse Training with Lipschitz Continuous Loss\\ Functions and a Weighted Group $l_0$-norm Constraint}

\author{Michael R. Metel\thanks{michael.metel@huawei.com}}

\affil{Huawei Noah's Ark Lab, Montr\'eal, Qc, Canada}

\maketitle	

\begin{abstract}%   <- trailing '%' for backward compatibility of .sty file
This paper is motivated by structured sparsity for deep neural network training. 
We study a weighted group $l_0$-norm constraint, and present the projection and normal cone of this set. Using randomized smoothing, we develop zeroth and first-order algorithms for minimizing a Lipschitz continuous function constrained by any closed set which can be projected onto. Non-asymptotic convergence guarantees are proven in expectation for the proposed algorithms for two related convergence criteria which can be considered as approximate stationary points. Two further methods are given using the proposed algorithms: one with non-asymptotic convergence guarantees in high probability, and the other with asymptotic guarantees to a stationary point almost surely. We believe in particular that these are the first such non-asymptotic convergence results for constrained Lipschitz continuous loss functions.
\end{abstract}

%\begin{keywords}
%  structured sparsity, L0-norm, Lipschitz continuity, stochastic optimization, non-asymptotic convergence 
%\end{keywords}

\section{Introduction}
\label{int}
This paper focuses on training deep neural networks with structured sparsity using a constrained optimization approach. A structured sparsity constraint allows for a simpler neural network architecture to be selected which can be deployed on low-resource devices. Research on sparsity in deep learning is vast, for a thorough background see \citep{hoefler2021}. Though much of this research is of a heuristic nature, our focus is on algorithms with theoretical convergence guarantees. The problem is modelled as  
\begin{alignat}{6}\label{eq:0}
	&&\min\limits_{w\in\RR^d}\text{ }&f(w)\quad&&\st\text{ }&w\in C
\end{alignat}
where $C$ is a weighted group $l_0$-norm constraint defined in Section \ref{constraint}. This work examines the Euclidean projection operator and the normal cone of $C$, and develops new non-asymptotic convergence results for general zeroth and first-order stochastic projected algorithms for assumptions on $f(w)$ applicable for a wide range of architectures in deep learning \citep{davis2020}. In particular, the function $f$ is only assumed to be Lipschitz continuous on a compact set, taking the form of the expected value of an integrable stochastic loss function $F(w,\xi)$, 
$f(w):=\EE[F(w,\xi)]$ , where $\xi\in\RR^p$ is a random vector from a probability space $(\Omega,\F,P)$. In the context of supervised learning, given samples $\xi^i=(x^i,y^i)$ for $i=1,2,...,\M$, where $\{x^i\}$ is a feature set, $\{y^i\}$ is a label set, and $F(w,\xi^i)$ is the loss associated with sample $i$, $f(w)$ can be replaced in \eqref{eq:0} by its approximation $\hat{f}(w)=\M^{-1}\sum_{i=1}^{\M} F(w,\xi^i)$.\\

The next section summarizes the required definitions and notation which will be used throughout the paper. Section \ref{constraint} presents the weighted group $l_0$-norm constraint. Section \ref{relworks} gives an overview of related works focusing on algorithms with theoretical convergence guarantees for Lipschitz continuous loss functions. Section \ref{smoothing} gives the detailed assumptions on $F(w,\xi)$, and presents the technique of using randomized smoothing to overcome the non-differentiability of the loss function. In Section \ref{constraintprop} the Euclidean projection operator and the normal cone for the proposed constraint set are given. Section \ref{trainalg} presents the Stochastic Projected Algorithm (SPA), which has a zeroth and a first-order version, with new non-asymptotic convergence results for two related convergence criteria, and a method using SPA which has an asymptotic convergence guarantee to a stationary point almost surely. Section \ref{backprop} shows how backpropagation can be used in conjunction with the first-order version of SPA for a wide range of deep learning architectures and validates its use in Section \ref{exper} where the theory of Section \ref{trainalg} is applied to train a neural network. Section \ref{conclusion} concludes the work. All proofs of results can be found in the Appendices A-E.

\section{Preliminaries}
\label{pre}  
For a set $S$, let the notation $x\xrightarrow{S} w$ mean $x\rightarrow w$ with $x\in S$, and for a discontinuous function $h$, $x\xrightarrow{h} w$ indicates that $x\rightarrow w$ with $h(x)\rightarrow h(w)$. For a function $h:\RR^d\rightarrow \RR$, 
\begin{alignat}{6}
	\limsup_{x\rightarrow w}h(x):=\inf_{\gamma>0}\bigg(\sup_{0<||x-w||_2<\gamma}h(x)\bigg).\nonumber
\end{alignat}
For a set-valued mapping $G: \RR^d\rightrightarrows \RR^d$,  
\begin{alignat}{6}
	&\Limsup_{x\rightarrow w}G(x):=\{y\in\RR^d:\exists \text{ sequences }x_k\rightarrow w \text{ and } y_k\rightarrow y \text{ with } y_k\in G(x_k) \text{ }\forall k\in \NN\}.\nonumber
\end{alignat}
For $w\in S$, the Fr\'echet normal cone equals
\begin{alignat}{6}
	\widehat{N}(w,S):=\{y\in \RR^n: \limsup_{x\xrightarrow{S} w} \frac{\langle y,x-w\rangle}{||x-w||_2}\leq 0\}\nonumber
\end{alignat}
and the Mordukhovich normal cone equals 
\begin{alignat}{6}
	N(w,S)=\Limsup_{x\xrightarrow{S}w}\widehat{N}(x,S).\nonumber
\end{alignat}
When $w\notin S$, $N(w,S)=\widehat{N}(w,S):=\{\emptyset\}$. For more information about normal cones see for example \cite{mordukhovich2013}.\\

For an extended real-valued function $h:\RR^d\rightarrow \RR\cup\{\pm\infty\}$, when finite, let $\widehat{\partial} h(w)$ denote its Fr\'echet subdifferential, defined as 
\begin{alignat}{6}
	&\widehat{\partial} h(w):=\{y\in\RR^n:\liminf_{x\rightarrow w} \frac{f(x)-f(w)-\langle y,x-w\rangle}{||x-w||_2}\geq 0\},\nonumber
\end{alignat}
and let $\partial h(w)$ denote its Mordukhovich subdifferential,
\begin{alignat}{6}
	\partial h(w):=\Limsup_{x\xrightarrow{h}w} \widehat{\partial}h(x).\nonumber
\end{alignat}
Assuming that $f(w)$ is locally Lipschitz continuous and $S$ is closed, a necessary condition for $\overline{w}$ to be locally optimal for the problem
\begin{alignat}{6}
	&&\min\limits_{w\in\RR^d}\text{ }&f(w)\quad&&\st\text{ }&w\in S\nonumber
\end{alignat}
is for \citep[Theorem 8.15 \& 9.13]{rockafellar2009}
\begin{alignat}{6}
	0\in \partial f(\overline{w})+N(\overline{w},S).\label{optcond}
\end{alignat}
As a non-asymptotic convergence criterion for optimization algorithms, for an $\epsilon>0$, $\overline{w}$ is an $\epsilon$-stationary point when 
\begin{alignat}{6}
	\dist(0,\partial f(\overline{w})+N(\overline{w},S))\leq \epsilon.\nonumber
\end{alignat}

Let $\overline{\partial}f(w)$ denote the Clarke subdifferential which equals $\overline{\partial}f(w)=\co\{\partial f(w)\}$, where $\co\{\cdot\}$ is the convex hull, given that $f(w)$ is locally Lipschitz continuous \citep[Theorem 9.61]{rockafellar2009}. When $f(w)$ is Clarke regular, meaning that its one-sided directional derivative exists and for all $v\in\RR^d$ $f'(w;v)=\max\limits_{g\in\overline{\partial} f(w)}\langle g,v\rangle$ \citep[Proposition 2.1.2 (b) \& Definition 2.3.4]{clarke1990}, the Clarke subdifferential coincides with the Fr\'echet and Mordukhovich subdifferentials \citep[Theorem 9.61 \& Corollary 8.11]{rockafellar2009}.\\

In this work, we will consider a relaxed version of \eqref{optcond}, which we call a Clarke-Mordukhovich (C-M) stationary point:
\begin{alignat}{6}
	0\in \overline{\partial} f(\overline{w})+N(\overline{w},S).\label{cmstat}
\end{alignat}

Let $B(x,r):=\{x+z:\norm{z}_2< r\}$ be the open Euclidean ball centered at $x$ with radius $r$, let $\overline{B}(x,r)$ be the corresponding closed Euclidean ball, and let $\overline{B}_{r}:=\{z\in\RR^d:|z_i|\leq r \text{ for }i=1,2,...,d\}$ denote the closed $l_{\infty}$-ball with radius $r>0$ centered at $0$.\\ 

We will also consider the Clarke $\epsilon$-subdifferential, 
\begin{alignat}{6}  
	\overline{\partial}_{\epsilon} f(w):=\text{co}\{\overline{\partial} f(x): x\in \overline{B}(w,\epsilon)\},\nonumber
\end{alignat}
which was introduced in \citep{goldstein1977}. This type of subdifferential has recently been used in the non-asymptotic convergence analysis of minimization algorithms for unconstrained Lipschitz continuous functions, see \citep{zhang2020,metel2021,kornowski2021} for more background. Besides its use for non-asymptotic convergence analysis, it holds that $\lim\limits_{\epsilon\rightarrow 0}\overline{\partial}_{\epsilon} f(\overline{w})=\overline{\partial} f(\overline{w})$ \citep[Lemma 7]{zhang2020}, which motivates the proposed C-M stationary point \eqref{cmstat} for our asymptotic convergence analysis.\\

The indicator function of a set $S$ equals 
\begin{alignat}{6}
	\delta_{S}(w)=\begin{cases}
		0 & \text{ if } w\in S\\ 
		\infty & \text{ otherwise,}\\ 
	\end{cases}\nonumber
\end{alignat}
and $2^{S}$ denotes its power set. For a random variable $X$, let $P_{X}$ denote the probability measure induced by the random variable $X$, i.e. for a Borel set $S$, $P_{X}(S)=P(\{\omega\in\Omega: X(\omega)\in S\})$. For an $n\in\NN$, let $[n]:=\{1,2,...,n\}$ and $[n]_{-1}:=\{0,1,...,n-1\}$. When studying the computational complexity of algorithms we will use the notation $\tilde{O}$ which is the standard big O notation with logarithmic terms ignored, e.g. $\log^k(x)=\tilde{O}(1)$ for any $k\in\RR$. 

\section{Weighted group $l_0$-norm constraint} 
\label{constraint}
The $l_0$-norm counts the number of non-zero elements in a vector $w\in\RR^d$, 
\begin{alignat}{6}
	||w||_0:=\sum_{i=1}^d\ind_{\{w\in \RR^d: w_i\neq 0\}}(w).\nonumber
\end{alignat} 
For an $n\leq d$, let $\{w^i\}_{i=1}^n$ be a partition of $w$, where $w^i$ is of dimension $d_i$ for each $i\in [n]$ and $\sum_{i=1}^n d_i=d$.
The weighted group $l_0$-norm constraint is then defined as 
\begin{alignat}{6}
	&C:=\{w\in\RR^d:\sum_{i=1}^np_i\ind_{\{w^i\neq 0\}}(w)\leq m\},\nonumber
\end{alignat}
where $p_i>0$ is a finite penalty associated with the subset of decision variables $w^i$, $\{w^i\neq 0\}$ denotes the set ${\{w\in \RR^d: \exists j\in[d_i],\text{ }w^i_j\neq 0\}}$, and $m>0$ is the maximum allowable aggregate penalty. The choice of the partition can be made to simplify a neural network's architecture, for example each $w^i$ can be the weights and bias of a neuron in a fully connected layer or of a filter in a convolutional layer. If $m$ is an upper bound on the available memory to store $w$ on a device, then each $p_i$ can be the amount of memory required for each $w^i$. We assume that each $p_i\leq m$. If there exists a $p_i>m$, the associated decision variables $w^i$ can be removed without affecting problem \eqref{eq:0}. For the Euclidean projection operator $\Pi_{C}(\cdot)$ to be nonempty, it is sufficient that $C$ is a closed set \citep[Example 1.20]{rockafellar2009}, which is verified in the next proposition.

\begin{proposition}
	\label{closedset}
	$C$ is a closed set.
\end{proposition}

\section{Related works}
\label{relworks}
The projection and normal cone of $C$ are presented in Section \ref{constraintprop}, which is an extension of the analysis of the $l_0$-norm constraint in \citep{bauschke2014}. Non-asymptotic convergence to an expected $\epsilon$-stationary point has been established for the proximal mini-batch SGD algorithm under the assumption that the function $f$ has a Lipschitz continuous gradient in \citep{xu2019non}. In general, neural networks are not differentiable so this result cannot be applied. More appropriate for deep learning optimization is the assumption that $f(w)$ is (locally) Lipschitz continuous.\\

For asymptotic convergence results, in \citep{davis2020}, the stochastic subgradient algorithm is proven to converge asymptotically to a Clarke stationary point almost surely for locally Lipschitz functions which admit a Whitney stratifiable graph, for step-sizes approaching zero in the limit, with an extension to the proximal stochastic subgradient algorithm. In \citep{bianchi2020}, the authors consider a fixed step-size and model the randomness of stochastic gradients in a manner more congruent with using SGD for locally Lipschitz loss functions, and prove a convergence result in probability to the set of Clarke stationary points. The authors also consider a projected SGD algorithm, in particular for closed Euclidean balls, which ameliorates some technical assumptions. 
A locally Lipschitz continuous generalized-differentiable \citep{norkin1980} function with a convex and closed constraint is considered in \citep{ruszczynski2020}. Asymptotic convergence to a Clarke stationary point for a stochastic subgradient method with averaging is proven.\\ 

Non-asymptotic convergence for a zeroth-order algorithm is presented in \citep{nesterov2017} for the minimization of deterministic Lipschitz continuous functions using Gaussian smoothing. Non-asymptotic convergence results for first-order methods in terms of the Clarke $\epsilon$-subdifferential, in the deterministic and stochastic setting are given in \citep{zhang2020}, under the assumption that loss functions are directionally-differentiable, and in the stochastic setting in \citep{metel2021} using iterate perturbation. A comparison of our convergence criteria and computational complexity is given in Section \ref{trainalg}.

\section{Randomized smoothing of $f(w)$}
\label{smoothing}
To overcome the non-differentiability of $f(w)$, the original problem can be replaced by a smoothed approximation (see Proposition \ref{gradcont}),
\label{randfunc}
\begin{alignat}{6}
	&\min\limits_{w\in \RR^d}&&\text{ }f_{\alpha}(w)\quad&\st &&\text{ }w\in C\cap \overline{B}_{\beta},\nonumber
\end{alignat}
where $f_{\alpha}(w):=\EE[f(w+u)]$ for a random vector $u:\Omega\rightarrow \RR^d$ uniformly distributed over 
$\overline{B}_{\frac{\alpha}{2}}$ for an $\alpha>0$. All $u_i$ are mutually independent random variables with marginal probability distributions equal to  
$$P_{u_i}=\begin{cases}
	\frac{1}{\alpha} & \text{ if } |u_i|\leq \frac{\alpha}{2}\\ 
	0 & \text{ otherwise.}  
\end{cases}\nonumber$$ 

The added constraint $\overline{B}_{\beta}$ for a $\beta>0$ is to allow us to assume that $f(w)$ is only Lipschitz continuous over a compact set around zero. If $f(w)$ is Lipschitz continuous over $\RR^d$, this constraint can be removed by setting $\beta=\infty$. The assumptions on $F(w,\xi)$ are similar to those used in \citep{metel2021}. For a $\kappa>\beta+\frac{\alpha}{2}$, we assume that $F(w,\xi)$ is a $\B_{\overline{B}_{\kappa}\times\RR^p}$-measurable function, where $\B_{(\cdot)}$ denotes the Borel $\sigma$-algebra. We assume that for each $\xi\in\RR^p$, $F(w,\xi)$ is continuous in $w\in \overline{B}_{\kappa}$, and for a measurable function $L_0(\xi)$, $F(w,\xi)$ is $L_0(\xi)$-Lipschitz continuous,
\begin{alignat}{6}
	|F(w,\xi)-F(w',\xi)|\leq L_0(\xi)||w-w'||_2,\label{loclip}
\end{alignat}
for all $w,w'\in \overline{B}_{\kappa}$ and for all $\xi\in\RR^p$ outside of a Borel null set. It is assumed that $L_0(\xi)$ is square integrable, $Q:=\EE[L_0(\xi)^2]<\infty$. It follows that $f(w)$ is Lipschitz continuous in $w\in \overline{B}_{\kappa}$.

\begin{proposition}
	\label{lipschitz_f}
	The function $f$ is  $L_0:=\EE[L_0(\xi)]$-Lipschitz continuous over $\overline{B}_{\kappa}$. 
\end{proposition}

Given that $f(w)$ is Lipschitz continuous over $w\in \overline{B}_{\kappa}$, it is differentiable almost everywhere over $w\in \overline{B}_{\beta+\alpha/2}$ by Rademacher's theorem \citep[Theorem 3.1]{heinonen2005}. The function $\nabla f$ may not be defined on a null set in $\overline{B}_{\beta+\alpha/2}$, so we define $\widetilde{\nabla}f(w)$ to be a $\B_{\overline{B}_{\beta+\alpha/2}}$-measurable function which for every $w\in \overline{B}_{\beta}$ equals 
$\nabla f(w+u)$ for almost every $u$ over $(\overline{B}_{\alpha/2},\B_{\overline{B}_{\alpha/2}},P_u)$. Similarly, the function $F(w,\xi)$ is differentiable almost everywhere over the product measure space 
${(\overline{B}_{\beta+\alpha/2}\times\RR^p,\B_{\overline{B}_{\beta+\alpha/2}\times\RR^{p}},m\times P_{\xi})}$, where $m$ is the Lebesgue measure restricted to Borel sets \citep[Property 1]{metel2021}. We define $\widetilde{\nabla}F(w,\xi)$ to be a $\B_{\overline{B}_{\beta+\alpha/2}\times\RR^p}$-measurable function, which for every $w\in \overline{B}_{\beta}$ equals $\nabla F(w+u,\xi)$ for almost every $(u,\xi)$ over $(\overline{B}_{\alpha/2}\times\RR^p,\B_{\overline{B}_{\alpha/2}\times\RR^{p}},P_u\times P_{\xi})$. Applying the results of \citep{bolte2021}, it is verified in Section \ref{backprop} that the output of backpropagation has the key properties of $\widetilde{\nabla}F(w,\xi)$, namely measurability and being equal to $\nabla F(w,\xi)$ almost everywhere for conditions which are widely applicable for deep learning applications.\\ 

Another approach to overcome the non-differentiability of $f(w)$ is to consider a zeroth-order algorithm. As proposed in \citep{gupal1977}, an unbiased stochastic estimation of the gradient of $f_{\alpha}(w)$ can be computed using the following finite-difference functions, $df: \RR^{2d}\rightarrow \RR^d$ and $dF: \RR^{2d+p}\rightarrow \RR^d$, defined 
component-wise as
\begin{alignat}{6}
	df_{i}(w,u_{\setminus i}):=&f(w_1+u_1,...,w_{i-1}+u_{i-1},w_i+\frac{\alpha}{2},w_{i+1}+u_{i+1},...,w_d+u_d)\nonumber\\
	-&f(w_1+u_1,...,w_{i-1}+u_{i-1},w_i-\frac{\alpha}{2},w_{i+1}+u_{i+1},...,w_{d}+u_d),\nonumber
\end{alignat}
and
\begin{alignat}{6}
	dF_{i}(w,u_{\setminus i},\xi):=&F(w_1+u_1,...,w_{i-1}+u_{i-1},w_i+\frac{\alpha}{2},w_{i+1}+u_{i+1},...,w_d+u_d,\xi)\nonumber\\
	-&F(w_1+u_1,...,w_{i-1}+u_{i-1},w_i-\frac{\alpha}{2},w_{i+1}+u_{i+1},...,w_{d}+u_d,\xi),\nonumber
\end{alignat}
where $u_{\setminus i}:=[u_1,...,u_{i-1},u_{i+1},...,u_d]^T$.\footnote{We use the notation $df(w,u)$ and $dF(w,u,\xi)$, but then switch to $df_{i}(w,u_{\setminus i})$ and $dF_{i}(w,u_{\setminus i},\xi)$ to make it clear that there is no $u_i$ argument for the $i^{th}$ component function.} 
\subsection{Properties of $f_{\alpha}(w)$}
Given the assumptions made about the use of randomized smoothing and the stochastic function $F(w,\xi)$, the following are resulting properties of $f_{\alpha}(w)$ which will be useful in Section \ref{trainalg} for the analysis of the proposed training algorithms. For some similar results when $f(w)$ is Lipschitz continuous over $\RR^d$ and the random vector $u$ is Gaussian, see \citep{nesterov2017}. Gupal (\citeyear{gupal1977}) motivated the use of uniform perturbation over $l_{\infty}$-balls, where similar results to Propositions \ref{grad} and \ref{gradcont} can be found.\\

The following proposition proves that unbiased estimates of $\nabla f_{\alpha}(w)$ can be generated using 
$\widetilde{\nabla}f$, $\widetilde{\nabla}F$, $df$, or $dF$ with samples of $u$ and $\xi$.
\begin{proposition}
	\label{grad}	
	For all $w\in \overline{B}_{\beta}$, 
	\begin{alignat}{6}
		\nabla f_{\alpha}(w)&=\EE[\widetilde{\nabla} f(w+u)]&=&\EE[\widetilde{\nabla} F(w+u,\xi)]\nonumber\\
		&=\alpha^{-1}\EE[df(w,u)]&=&\alpha^{-1}\EE[dF(w,u,\xi)].\nonumber
	\end{alignat}	
\end{proposition}
The next proposition relates the gradient of $f_{\alpha}(w)$ with the Clarke $\epsilon$-subdifferential of $f(w)$. 
\begin{proposition}
	\label{clarkeps}	
	Assume that $\kappa> \beta + \sqrt{d}\frac{\alpha}{2}$. For all $w\in \overline{B}_{\beta}$ with $\widehat{\alpha}=\sqrt{d}\frac{\alpha}{2}$, 
	$\nabla f_{\alpha}(w)\in \overline{\partial}_{\widehat{\alpha}} f(w)$.	
\end{proposition}
Proposition \ref{clarkeps} required a stronger condition on $\kappa$ to ensure that  $\overline{\partial}_{\widehat{\alpha}} f(w)$ is well-defined, meaning that the Clarke subdifferential is only being considered for values of $x\in \RR^d$ where $f(w)$ is Lipschitz continuous on a neighbourhood of $x$. The following proposition focuses on properties of $f_{\alpha}(w)$ and its relation to $f(w)$. 

\vbox{\begin{proposition}
		\label{lipcont}
		\text{ }
		\begin{enumerate}
			\item $f_{\alpha}(w)$ is $L_0$-Lipschitz continuous for $w\in \overline{B}_{\beta}$.
			\item For all $w\in \overline{B}_{\beta}$,
			$|f_{\alpha}(w)-f(w)|\leq \alpha L_0\sqrt{\frac{d}{12}}$.	
			\item For a closed set $S$, let $w_{\alpha}^*$ and $w^*$ be minimizers of $f_{\alpha}(w)$ and $f(w)$ respectively for $w\in S\cap \overline{B}_{\beta}$, then 
			$|f_{\alpha}(w_{\alpha}^*)-f(w^*)|\leq \alpha L_0\sqrt{\frac{d}{12}}$.
			\item For any two values $w,w'\in \overline{B}_{\beta}$,
			$|f_{\alpha}(w)-f_{\alpha}(w')|\leq 2\beta\sqrt{d}L_0$.
		\end{enumerate}
\end{proposition}}
The following proposition gives the Lipschitz constant of $\nabla f_{\alpha}(w)$, and will be referred to as the smoothness of $f_{\alpha}(w)$.
\begin{proposition}
	\label{gradcont}
	For all $w\in \overline{B}_{\beta}$,
	$\nabla f_{\alpha}(w)$ is $2\alpha^{-1}\sqrt{d}L_0$-Lipschitz continuous.
\end{proposition}

Considering the sample mean of a mini-batch of estimators of $\nabla f_{\alpha}(w)$, the next proposition gives bounds on the trace of their covariance matrices and on the expected value of their squared $l_2$-norm, which will be used in the convergence analysis of Section \ref{trainalg}.
\begin{proposition}
	\label{dq}
	For all $w\in \overline{B}_{\beta}$,
	\begin{enumerate}
		\item $\EE[||\nabla f_{\alpha}(w)-\frac{1}{M\alpha}\sum_{i=1}^MdF(w,u^i,\xi^i)||^2_2]\leq \frac{dQ}{M}$
		\item$\EE[||\nabla f_{\alpha}(w)-\frac{1}{M}\sum_{i=1}^M\widetilde{\nabla} F(w+u^i,\xi^i)||^2_2]\leq \frac{Q}{M}$
		\item$\EE[||\frac{1}{M\alpha}\sum_{i=1}^MdF(w,u^i,\xi^i)||^2_2]\leq dQ$
		\item$\EE[||\frac{1}{M}\sum_{i=1}^M\widetilde{\nabla} F(w+u^i,\xi^i)||^2_2]\leq Q$,
	\end{enumerate}
	where $\{u^i\}$ and $\{\xi^i\}$ are independent samples of $u$ and $\xi$.
\end{proposition}

\section{Properties of $C\cap \overline{B}_{\beta}$}
\label{constraintprop}
In this section we give the projection onto $C\cap \overline{B}_{\beta}$, and its Fr\'echet and Mordukhovich normal cones. For the projection and normal cone of the set $\{w\in\RR^d:||w||_0\leq m\}$, see \citep{bauschke2014}.
The projection onto $C\cap \overline{B}_{\beta}$ requires solving a 0-1 knapsack problem. This problem is NP-complete, though it can be solved in pseudo-polynomial time when all $p_i\in\ZZ_{>0}$ and $m\in\ZZ_{>0}$, which holds when allocating memory as described in Section \ref{constraint}. For further background on this problem and algorithms see \citep{kellerer2004}. 

\subsection{Projection onto $C\cap \overline{B}_{\beta}$}

The next proposition shows how the projection onto $C\cap \overline{B}_{\beta}$ can be computed using a 0-1 knapsack problem.

\begin{proposition}
	\label{projection}
	For any $w\in\RR^d$, let $Z^*$ equal the set of optimal solutions of the 0-1 knapsack problem, 
	\begin{alignat}{6}
		\max_{z\in\{0,1\}^n}\text{ }&\sum_{i=1}^nz_i(||w^i||^2_2-||\max(|w^i|-\beta,0)||^2_2)\label{knap2}\\
		\text{s.t. }&\sum_{i=1}^nz_ip_i\leq m,\nonumber
	\end{alignat}
	where $||\max(|w^i|-\beta,0)||^2_2:=\sum_{j=1}^{d_i}(\max(|w^i_j|-\beta,0))^2$. The projection $\Pi_{C\cap \overline{B}_{\beta}}(w)$ equals 
	\begin{alignat}{6}
		\Pi_{C\cap \overline{B}_{\beta}}(w)=\{x\in\RR^d: \exists z^*\in Z^*,\text{ } 
		&x^i=\sgn(w^i)\min(|w^i|,\beta)\text{ if } z^*_i=1,\nonumber\\
		&x^i=0 \text{ otherwise }\forall i\in [n]\},\nonumber
	\end{alignat}
	where $x^i=\sgn(w^i)\min(|w^i|,\beta)$ denotes ${x_j^i=\sgn(w_j^i)\min(|w^i_j|,\beta)}$ for $j=1,2,...,d_i$. 
\end{proposition}

The following remark gives the projection onto $C$, which can be verified by taking $\beta\rightarrow \infty$ in Proposition \ref{projection}. 

\begin{remark}
	\label{procases}
	For the projection onto $C$ for any $w\in\RR^d$, the 0-1 knapsack problem \eqref{knap2} becomes  
	\begin{alignat}{6}
		\max_{z\in\{0,1\}^n}\text{ }&\sum_{i=1}^nz_i||w^i||^2_2\label{eq:26}\\
		\text{s.t. }&\sum_{i=1}^nz_ip_i\leq m.\nonumber
	\end{alignat}
	If $Z^*$ equals the set of optimal solutions of \eqref{eq:26}, then the projection $\Pi_C(w)$ equals
	\begin{alignat}{6}
		\Pi_{C}(w)=\{&x\in\RR^d: \exists z^*\in Z^*,\text{ }
		x^i=w^i\text{ if } z^*_i=1,\text{ }x^i=0 \text{ otherwise }\forall i\in [n]\}.\nonumber
	\end{alignat} 	
\end{remark}

\subsection{Normal cones of $C\cap \overline{B}_{\beta}$}
Assume that $w\in C\cap \overline{B}_{\beta}$, and let $I(w):=\{i\in [n]:w^i\neq 0\}$ be the indices of the subsets of non-zero weights, and let $J(w):=\{j\in [n]\setminus I(w): \sum_{i\in I(w)}p_i +p_j\leq m\}$ be the indices of subsets which are zero, but are not constrained to be. The following proposition gives the Fr\'echet normal cone to the set $C\cap \overline{B}_{\beta}$.
    
\begin{proposition}
	\label{frechnorm}
	For any $w\in C\cap \overline{B}_{\beta}$, 	
	\begin{alignat}{6}
		\widehat{N}(w,C\cap \overline{B}_{\beta})=\bigg\{&y \in \RR^d : \forall i\in I(w)\cup J(w), \forall j\in[d_i],\text{ }
		&y^i_j\in
		\begin{cases}
			\RR_{\geq 0}&\text{if }w^i_j=\beta\\
			0&\text{if }|w^i_j|<\beta\\
			\RR_{\leq 0}&\text{if }w^i_j=-\beta
		\end{cases}\bigg\}.\label{freeq}		
	\end{alignat} 
\end{proposition}

The following remark gives the Fr\'echet normal cone to the set $C$, which can be verified by taking $\beta\rightarrow \infty$ in Proposition \ref{frechnorm}. 

\begin{remark}	
	For any $w\in C$, 	
	\begin{alignat}{6}
		\widehat{N}(w,C)=\{y \in \RR^d : \forall i\in I(w)\cup J(w),\text{ }y^i=0\}.\nonumber
	\end{alignat}	
\end{remark}

Let $Y:=\{X\subseteq2^{[n]}: I(w)\subseteq X,\text{ }\sum_{i\in X}p_i\leq m\text{ and } \sum_{i\in X}p_i+p_j> m\text{ }\forall j\notin X\}$, which contains all of the sets of indices $X$ containing $I(w)$ which make the constraint $\sum_{i\in X}p_i\leq m$ tight, in the sense that no further feasible index can be added to $X$. The next proposition gives the Mordukhovich normal cone to the set $C\cap \overline{B}_{\beta}$.

\begin{proposition}
	\label{limitnorm}	
	For any $w\in C\cap \overline{B}_{\beta}$, 	
	\begin{alignat}{6}
		N(w,C\cap \overline{B}_{\beta})=\bigg\{&y \in \RR^d : \exists X\in Y, \forall i\in X, \forall j\in[d_i],\text{ }
		&y^i_j\in
		\begin{cases}
			\RR_{\geq 0}&\text{if }w^i_j=\beta\\
			0&\text{if }|w^i_j|<\beta\\
			\RR_{\leq 0}&\text{if }w^i_j=-\beta
		\end{cases}\bigg\}.\label{moreq}
	\end{alignat}	
\end{proposition}

The next remark gives the Mordukhovich normal cone to the set $C$. 

\begin{remark}
	For any $w\in C$, 
	\begin{alignat}{6}
		N(w,C)=\{y \in \RR^d : \exists X\in Y,\text{ }y^i=0\text{ }\forall i\in X\}.\nonumber
	\end{alignat}	
\end{remark}

\section{Training Algorithm}
\label{trainalg}

\begin{algorithm}[H]
	\caption{Stochastic Projected Algorithm (SPA)} 
	\begin{algorithmic}
		\STATE {\bfseries Input:} $w^1\in S\cap \overline{B}_{\beta}$, $\eta>0$, $K\in\ZZ_{>0}$,  $M\in\ZZ_{>0}$	
		\STATE $R\sim\text{uniform}\{2,...,K+1\}$		
		\FOR{$k=1,2,...,R-1$} 
		\STATE Sample $u^{k,i}\sim P_u$ for $i=1,...,M$
		\STATE Sample $\xi^{k,i}\sim P_{\xi}$ for $i=1,...,M$
		\STATE (1) ${w^{k+1}\in\Pi_{S\cap \overline{B}_{\beta}}(w^k-\frac{\eta}{M\alpha}\sum_{i=1}^MdF(w^k,u^{k,i},\xi^{k,i}))}$
		\STATE OR
		\STATE (2) ${w^{k+1}\in\Pi_{S\cap \overline{B}_{\beta}}(w^k-\frac{\eta}{M}\sum_{i=1}^M\widetilde{\nabla} F(w^k+u^{k,i},\xi^{k,i}))}$		
		\ENDFOR
		\STATE {\bfseries Output:} $w^R$ 
	\end{algorithmic}
	\label{alg1}
\end{algorithm}
The following convergence results of SPA (Algorithm \ref{alg1}) are applicable for any constraint set $S$ which is closed and for which there exists a computable element of the Euclidean projection onto $S\cap \overline{B}_{\beta}$. If \eqref{loclip} holds for all $w,w'\in\RR^d$, then the projection operator can be simplified to $\Pi_{S}$. The proof of Theorem \ref{maincon} is an adaptation of the proof of \citep[Theorem 2]{xu2019non}. For some similar results of this section for a first-order algorithm for unconstrained problems see \citep{metel2021}. SPA has two settings: (1) is a zeroth-order and (2) is a first-order algorithm. The algorithm requires that $w^1\in S\cap \overline{B}_{\beta}$. A simple choice is to pick an arbitrary $w^0\in\RR^d$ and to set $w^1\in\Pi_{S\cap \overline{B}_{\beta}}(w_0)$, but when training a neural network with $S=C$ and a high sparsity level, there is a risk of initializing the neural network with layer collapse \citep[Page 20]{hoefler2021}, where all weights in a layer are set to zero, disconnecting the network. We highlight that the initial $w^1\in S\cap \overline{B}_{\beta}$ can be chosen to ensure that there are non-zero weights in each layer, or any other desired property.\\

We consider two convergence criteria. A solution $\overline{w}$ is an expected $(\epsilon_1,\epsilon_2)$-stationary point if 
\begin{alignat}{6}
	&|f_{\alpha}(w)-f(w)|\leq\epsilon_1\text{ for all }w\in \overline{B}_{\beta},\quad\text{and}\quad
	&\EE[\dist(0,\nabla f_{\alpha}(\overline{w})+N(\overline{w},S\cap \overline{B}_{\beta}))]\leq\epsilon_2,\nonumber
\end{alignat}
which guarantees that $\overline{w}$ is an expected $\epsilon_2$-stationary point for a smooth approximation of $f(w)$ with a uniform error from $f(w)$ within $\epsilon_1$. A solution $\overline{w}$ is an expected $(\epsilon_3,\epsilon_4)$-stationary point if 
\begin{alignat}{6}
	\widehat{\alpha}\leq\epsilon_3\quad\text{and}\quad
	\EE[\dist(0,\overline{\partial} f_{\widehat{\alpha}}(\overline{w})+N(\overline{w},S\cap \overline{B}_{\beta}))]\leq\epsilon_4, \nonumber
\end{alignat}
which can be seen as a relaxation of an expected $\epsilon_4$-stationary point, replacing the Mordukhovich with a Clarke $\widehat{\alpha}$-subdifferential. It will also be used in Theorem \ref{asymptotic} for a method with an asymptotic convergence guarantee to a C-M stationary point. These convergence criteria are related as for sufficiently small $\alpha$ an $(\epsilon_1,\epsilon_2)$-stationary point implies an $(\epsilon_3,\epsilon_4)$-stationary point with $\epsilon_3=\sqrt{d}\frac{\alpha}{2}$ and $\epsilon_4=\epsilon_2$ using Proposition \ref{clarkeps}.\\ 

The next proposition verifies the existence of a Borel measurable selection of the projection operator $\Pi_{S\cap \overline{B}_{\beta}}(\cdot)$, and of the measurability of the distance functions used in the convergence criteria. By the assumptions that $F(w,\xi)$ and $\widetilde{\nabla} F(w,\xi)$ are Borel measurable, the iterates $\{w^{k}\}$ from SPA are measurable using such a selection of $\Pi_{S\cap \overline{B}_{\beta}}(\cdot)$. This proposition also covers C-M stationary points, i.e. $\overline{\partial}_0 f(w)=\overline{\partial} f(w)$ \citep[Corollary 2.5]{goldstein1977}.

\begin{proposition}
	\label{meas}
	Let $S$ be a closed set. There exists a measurable selection of $\Pi_{S\cap \overline{B}_{\beta}}(\cdot)$. For any $0\leq\epsilon\leq\frac{\alpha}{2}$
	$\dist(0,\overline{\partial}_{\epsilon}f(w)+N(w,S\cap \overline{B}_{\beta}))$, and $\dist(0,\nabla f_{\alpha}(w)+N(w,S\cap \overline{B}_{\beta}))$ are Borel measurable functions in $w\in \overline{B}_{\beta}$.
\end{proposition} 

Together with Proposition \ref{lipcont}.2, the following theorem presents the non-asymptotic convergence to an expected $(\epsilon_1,\epsilon_2)$-stationary point of SPA.

\begin{theorem}	
	\label{maincon}	
	Let $\eta=\frac{\alpha}{3\rho\sqrt{d}L_0}$ for $\rho>0$ and let 
	$\tau\geq 0$ such that $\rho+\tau>1$. For a solution from SPA given any choice of $w^1\in S\cap \overline{B}_{\beta}$, $K\in\ZZ_{>0}$, and $M\in\ZZ_{>0}$, it holds that 
	\begin{alignat}{6}
		&\EE[\dist(0,\nabla f_{\alpha}(w^{R})+N(w^{R},S\cap \overline{B}_{\beta}))^2]\leq& C_1\frac{2\alpha^{-1}\sqrt{d} L_0\Delta}{K}+C_2\frac{\upsilon Q}{M},\label{theorem7}
	\end{alignat}
	where 
	\begin{alignat}{6}
		&&\Delta&\geq f_{\alpha}(w^1)-f_{\alpha}(w^*_{\alpha})\nonumber\\
		&&&\leq \min\bigg(2\beta\sqrt{d}L_0,f(w^1)-f(w^*)+\alpha L_0\sqrt{\frac{d}{3}}\bigg)\label{delbound}
	\end{alignat}
	for minimizers $w_{\alpha}^*$ and $w^*$ of $f_{\alpha}(w)$ and $f(w)$ respectively for $w\in S\cap \overline{B}_{\beta}$, $C_1:=\frac{2(1+3\rho)}
	{(\tau+\rho-1)}+3\rho$, ${C_2:=\frac{4(1+3\rho)}
		{(\tau+\rho-1)}\bigg(\frac{1}{2}+\frac{2}{3}\frac{M\tau}{\rho^2}\bigg)+3}$, and 
	${\upsilon:=\begin{cases}
			d & \text{ if using } (1)\\ 
			1 & \text{ if using } (2).\\ 
	\end{cases}}$
\end{theorem}
Inequality \eqref{delbound} gives valid choices for $\Delta$ which are easier to compute and related to the true loss function $f$. Using Theorem \ref{maincon} and Proposition \ref{clarkeps}, the following corollary holds.
\begin{corollary}	
	\label{mainconeps}	
	Assume that $\kappa> \beta + \sqrt{d}\frac{\alpha}{2}$ and $\widehat{\alpha}=\sqrt{d}\frac{\alpha}{2}$. Let $\eta$, $\rho$, $\tau$, $\Delta$, $C_1$, and $C_2$ be defined as in Theorem \ref{maincon}. For a solution from SPA given any choice of $w^1\in S\cap \overline{B}_{\beta}$, $K\in\ZZ_{>0}$, and $M\in\ZZ_{>0}$, it holds that 
	\begin{alignat}{6}
		&\EE[\dist(0,\overline{\partial} f_{\widehat{\alpha}}(w^{R})+N(w^{R},S\cap \overline{B}_{\beta}))^2]\leq C_1\frac{2\alpha^{-1}\sqrt{d} L_0\Delta}{K}+C_2\frac{\upsilon Q}{M}.\label{coroll8}
	\end{alignat}	
\end{corollary}
The parameter $\tau$ is not required and in fact $\tau>0$ results in the constant term 
\begin{alignat}{6}
	&&&\frac{4(1+3\rho)}
	{(\tau+\rho-1)}\bigg(\frac{2}{3}\frac{\tau}{\rho^2}\bigg)\upsilon Q\nonumber
\end{alignat}
in the expansion of the right-hand-side of \eqref{theorem7} and \eqref{coroll8}. The parameter $\tau$ is included so that the convergence bounds are applicable for any step-size $\eta>0$, though it will likely be poor unless $Q\approx 0$. A large $\sqrt{d}$, such as for deep neural networks, will result in a small step-size $\eta$ as well as require a large $K$ to get an adequate convergence guarantee. The inclusion of $\tau$ is also an attempt to remedy this when $Q\approx 0$, i.e. replacing $\rho>0$ and $\tau=0$ with $\rho'>0$ and $\tau'>0$ such that $\rho'+\tau'=\rho$ will increase $\eta$ and decrease $C_1$. For the remainder of this section we will assume that $\tau=0$, which results in $C_1$ and $C_2$ being equal to
$C_1^{\tau=0}:=\frac{2+3\rho+3\rho^2}{\rho-1}$ and $C_2^{\tau=0}:=\frac{9\rho-1}{\rho-1}$. Table \ref{t:1} presents some choices for $\rho$ resulting in $C_1^{\tau=0}$ and $C_2^{\tau=0}$ being integer-valued.
\begin{table}[H]
	\caption{Some choices for $\rho$.}
	\label{t:1}
	\begin{center}
		\tabulinesep=0.75mm
		\begin{tabu}{c|cccc}
			\hline
			%	\toprule
			$\rho$ & $\frac{4}{3}$&$\frac{5}{3}$&2&5\\
			\hline	
			%	\midrule
			$C_1^{\tau=0}$ & 34 &23 &20 & 23\\
			$C_2^{\tau=0}$ & 33 & 21 & 17 & 11 \\
			%	\bottomrule
			\hline
		\end{tabu}
	\end{center}
	\vspace{-0.5cm}
\end{table}

The following corollary gives the computational complexity to guarantee an expected $(\epsilon_1,\epsilon_2)$ or $(\epsilon_3,\epsilon_4)$-stationary point in terms of the number of either $dF(w^k,u^{k,i},\xi^{k,i})$ or $\widetilde{\nabla} F(w^k+u^{k,i},\xi^{k,i})$ computations, which will be referred to as {\it gradient calls}, and in terms of the number of projections.
\begin{corollary}	
	\label{complexity}	
	Running SPA as described in Theorem \ref{maincon} with  
	$\alpha=\frac{\epsilon_1}{L_0\sqrt{\frac{d}{12}}}$, $\tau=0$, 
	\begin{alignat}{6}
		&&K=\left\lceil C_1\sqrt{\frac{4}{3}}\frac{d L_0^2\Delta}{\epsilon_1\epsilon^2_2}\right\rceil,\quad\text{and}\quad M=\left\lceil C_2\frac{2\upsilon Q}{\epsilon_2^2}\right\rceil\nonumber
	\end{alignat}
	guarantees an expected $(\epsilon_1,\epsilon_2)$-stationary point. Assuming that $\kappa>\beta+\epsilon_3$, and setting 
	$\alpha=\frac{2\epsilon_3}{\sqrt{d}}$, $\tau=0$,
	\begin{alignat}{6}
		&&K=\left\lceil C_1\frac{2dL_0\Delta}{\epsilon_3\epsilon^2_4}\right\rceil,\quad\text{and}\quad M=\left\lceil C_2\frac{2\upsilon Q}{\epsilon_4^2}\right\rceil\nonumber
	\end{alignat}
	guarantees an expected $(\epsilon_3,\epsilon_4)$-stationary point. Using these choices of $\alpha$, $\tau$, $K$, and $M$ give gradient call complexities of $O(\epsilon_1^{-1}\epsilon_2^{-4})$ and $O(\epsilon_3^{-1}\epsilon_4^{-4})$, and projection operator complexities of $O(\epsilon_1^{-1}\epsilon_2^{-2})$ and $O(\epsilon_3^{-1}\epsilon_4^{-2})$ to achieve an expected $(\epsilon_1,\epsilon_2)$ and $(\epsilon_3,\epsilon_4)$-stationary point, respectively.
\end{corollary}

The next corollary gives the computational complexity for an $(\epsilon_1,\epsilon_2)$ or $(\epsilon_3,\epsilon_4)$-stationary point with a probability of at least $1-\gamma$ for any $\gamma\in(0,1)$ using the method proposed in \citep[Section 2.2]{ghadimi2013}. SPA is required to be run $r\in\ZZ_{>0}$ times, generating $r$ different solutions, with the result holding for the solution which minimizes the distance to stationarity using a sample mean approximation of $\nabla f_{\alpha}(w)$.
\begin{corollary}
	\label{comp2}	
	For any $\gamma\in(0,1)$ and $\epsilon_1,\epsilon_2>0$ or $\epsilon_3,\epsilon_4>0$, with $\kappa>\beta+\epsilon_3$, assume that SPA is run $r:=\lceil-\ln(c\gamma)\rceil$ times for any $c\in(0,1)$ according to Theorem \ref{maincon} with $\tau=0$,
	$\alpha=\frac{\epsilon_1}{L_0\sqrt{\frac{d}{12}}}$ or $\alpha=\frac{2\epsilon_3}{\sqrt{d}}$,
	\begin{alignat}{6}
		&K=\left\lceil C_1\sqrt{\frac{4}{3}}\frac{d L_0^2\Delta}{\epsilon_1(\epsilon_2')^2}\right\rceil\quad&&\text{and}\quad &&M=\left\lceil C_2\frac{2\upsilon Q}{(\epsilon_2')^2}\right\rceil,\quad\text{or}\nonumber\\
		&K=\left\lceil C_1\frac{2dL_0\Delta}{\epsilon_3(\epsilon'_4)^2}\right\rceil\quad&&\text{and}\quad &&M=\left\lceil C_2\frac{2\upsilon Q}{(\epsilon'_4)^2}\right\rceil,\nonumber
	\end{alignat}
	where $\epsilon_2'=\sqrt{\frac{\epsilon^2_2-6\psi\frac{Q}{T}}{4e}}$, $\epsilon_4'=\sqrt{\frac{\epsilon^2_4-6\psi\frac{Q}{T}}{4e}}$, $\psi=\frac{\lceil-\ln(c\gamma)\rceil}{(1-c)\gamma}$, $e:=\exp(1)$, and $T=\lceil6\phi\psi\frac{Q}{\epsilon^2_2}\rceil$ or $T=\lceil6\phi\psi\frac{Q}{\epsilon^2_4}\rceil$ for any $\phi>1$, outputting solutions $W:=\{w^1,...,w^{r}\}$. 
	Let $\{u^i\}_{i=1}^T$ and $\{\xi^i\}_{i=1}^T$ be independent samples of $u$ and $\xi$, and let $w^*\in W$ be chosen such that 
	\begin{alignat}{6}
		w^*\in\argmin\limits_{w\in W}\dist(0,G(w)+N(w,S\cap \overline{B}_{\beta})),\label{choosew}
	\end{alignat}  
	where 
	\begin{alignat}{6}
		G(w):=\begin{cases}
			\frac{1}{T\alpha}\sum_{i=1}^TdF(w,u^i,\xi^i) & \text{ if using } (1)\\ 
			\frac{1}{T}\sum_{i=1}^T\widetilde{\nabla} F(w+u^i,\xi^i) & \text{ if using } (2)\\ 
		\end{cases}\nonumber
	\end{alignat}	
	in SPA. It follows that $w^*$ is an $(\epsilon_1,\epsilon_2)$ or $(\epsilon_3,\epsilon_4)$-stationary point with a probability of at least $1-\gamma$,  
	it is generated with $\tilde{O}\left(\epsilon^{-1}_1\epsilon^{-4}_2+\gamma^{-1}\epsilon^{-2}_2\right)$ or
	$\tilde{O}\left(\epsilon^{-1}_3\epsilon^{-4}_4+\gamma^{-1}\epsilon^{-2}_4\right)$ gradient calls, and $\tilde{O}\left(\epsilon^{-1}_1\epsilon^{-2}_2\right)$ or
	$\tilde{O}\left(\epsilon^{-1}_3\epsilon^{-2}_4\right)$ projections.
\end{corollary}

The optimization problem \eqref{choosew} requires knowledge of the normal cone of $S\cap \overline{B}_{\beta}$ as given in Proposition \ref{limitnorm} for $S=C$, and is solved by computing the distance $\dist(0,G(w)+N(w,S\cap \overline{B}_{\beta}))$ $\tilde{O}(1)$ times, once for each $w\in W$. The binary integer program discussed in the next proposition can be found in the proof, see \eqref{fin_binp}, but its requirement is only for pathological cases in neural network training when a solution equals $w^i=0$ for an $i\in[n]$ which is not constrained to be zero.
\begin{proposition}	
	\label{normcomp}
	If $Y=I(w)$ as defined above Proposition \ref{limitnorm}, $\dist(0,G(w)+N(w,C\cap \overline{B}_{\beta}))=||G(w)+v||_2$ where 
\begin{alignat}{6}
	v^i_j=\begin{cases}
		0 & \text{ if } i\in I(w) \text{ and } \neg U^i_j\\
		-G^i_j(w) & \text{otherwise,}\\
	\end{cases}\nonumber
\end{alignat}
where $U^i_j:= (|w^i_j|=\beta) \land (\sgn(G^i_j(w))=-\sgn(w^i_j))$ for $i\in[n]$ and $j\in d_i$.
When there exists an $X\in Y$ such that $X\setminus I(w)\neq\{\emptyset\}$, assume that $\{p_i\}\subset\QQ_{>0}$. The distance $\dist(0,G(w)+N(w,C\cap \overline{B}_{\beta}))$ can be computed by solving a binary integer program with $|[n]\setminus I(w)|$ binary variables.
\end{proposition}

\noindent{\bf Comparison of computational complexity}\\
Our gradient call complexity matches that of \citep[Section 7]{nesterov2017} to achieve an expected $(\epsilon_1,\epsilon_2)$-stationary point for an unconstrained deterministic function $f$. Their random gradient-free oracle only requires two function evaluations, whereas $df(w,u)$ requires $2d$ function evaluations. The Gaussian smoothing is computationally appealing but we would need to assume that $f(w)$ is Lipschitz continuous over $\RR^d$ as the function calls within $df(w,u)$ would now be evaluated at any point in $\RR^d$ given the expanded image of normal random variables compared to $u\in\overline{B}_{\frac{\alpha}{2}}$. In \citep{metel2021} the same computational complexity is proven for an expected $(\epsilon_3,\epsilon_4)$-stationary point for unconstrained stochastic functions. In \citep{zhang2020}, a better computational complexity of $\tilde{O}(\epsilon_3^{-1}\epsilon_4^{-3})$ is proven to achieve an $(\epsilon_3,\epsilon_4)$-stationary point in high probability in the deterministic setting. In the stochastic setting, \citep{zhang2020} proves a computational complexity of $\tilde{O}(\epsilon_3^{-1}\epsilon_4^{-4})$ to achieve an expected $(\epsilon_3,\epsilon_4)$-stationary point similar to our work.\\

The next theorem proves that the set of solutions from running SPA with increasing accuracy has an asymptotic convergence guarantee to a C-M stationary point almost surely.
\begin{theorem}	
	\label{asymptotic}
	Let $\{\epsilon^i_3\}$ and $\{\epsilon^i_4\}$ be strictly decreasing positive sequences approaching $0$ in the limit, with $\alpha^i$, $K^i$ and $M^i$ set to guarantee an expected $(\epsilon^i_3,\epsilon^i_4)$-stationary point running SPA according to Corollary \ref{complexity} assuming that $\kappa>\beta +\epsilon^1_3$. 
	Assume that SPA is run according to Theorem \ref{maincon} with $\tau=0$, $\alpha=\alpha^i$, $K=K^i$, and $M=M^i$ for $i=1,2,...$, giving solutions $\{w^i\}$. If $\beta$ is finite, there exists an accumulation point $\overline{w}$ of $\{w^i\}$ and it is a C-M stationary point almost surely.	Otherwise, any accumulation point $\overline{w}$ of $\{w^i\}$ is a C-M stationary point almost surely.
\end{theorem}

\section{Using Backpropagation}
\label{backprop}
This section considers computing  $\widetilde{\nabla}F(w,\xi)$ using backpropagation for a problem setting entailing a wide range of deep learning applications. This is demonstrated using the results of \citep{bolte2021}. Assume that $\xi$ maps to a countable number of values $\{\xi_i\}_{i=1}^{\infty}$ almost surely, $\PP(\xi\in \{\xi_i\}_{i=1}^{\infty})=1$, and assume that for each $\xi_k\in \{\xi_i\}_{i=1}^{\infty}$, $F(w,\xi_k)$ for $w\in \RR^d$ can be written as a composition of locally Lipschitz continuous functions $\{\sigma_{j}\}_{j\in I_k}$, for an index set $I_k$, and assume that the functions $\{\sigma_{j}\}_{j\in I_k}$ are definable in the same o-minimal structure. 

\begin{proposition}\citep[Corollary 5]{bolte2021}
	\label{boltebackprop}
	For each $\xi_k\in\{\xi_i\}_{i=1}^{\infty}$ set $\widetilde{\nabla}F(w,\xi_k)$ equal to the output of backpropagation using a measurable selection $\widetilde{\nabla}\sigma_{j}(\cdot)\in \overline{\partial}\sigma_{j}(\cdot)$, which exists, for all $j\in I_k$, and for $\xi\notin\{\xi_i\}_{i=1}^{\infty}$ set $\widetilde{\nabla}F(w,\xi)=a$ for any $a\in\RR^d$. $\widetilde{\nabla}F(w,\xi)$ is Borel measurable for $(w,\xi)\in\RR^{d+p}$ and equals the gradient of $F(w,\xi)$ for almost every $(w,\xi)\in\RR^{d+p}$.
\end{proposition} 

We focus on the o-minimal structure of the ordered real exponential field, 

\noindent$\RR_{\text{exp},<}:=\RR(+,\cdot,0,1,<,\exp)$, which provides a wide class of definable functions typically found in deep learning architectures. For a short background on o-minimal structures see for example \citep{wilkie2007}. The next proposition verifies the validity of using backpropagation for the building blocks used in the neural network considered in the next section, and also contains a sufficient background on o-minimal structures to understand the result. Other activation functions typically used in deep learning can be shown to have the following properties as well. We refer to what is computed during backpropagation as a {\it bp} gradient. Conv2d and MaxPool2d are defined as tensor-valued functions, but it is sufficient to consider their component functions separately. 
\begin{proposition}
	\label{omin}
	The affine map, ReLU, the component functions of Conv2d and MaxPool2d, and the loss function CrossEntropyLoss are definable in the o-minimal structure of $\RR_{\text{exp},<}$, and their bp gradients are measurable selections of their Clarke subdifferentials.
\end{proposition}

\section{Training a Neural Network}
\label{exper}

We trained a Lenet-5 type neural network on the MNIST (MN) and FashionMNIST (FMN) datasets constrained by $C\cap \overline{B}_{\beta}$. The projection operator $\Pi_{C\cap \overline{B}_{\beta}}(\cdot)$ was computed using a branch-and-bound (BNB) algorithm to solve the $0-1$ knapsack problem \eqref{knap2}. A sampling approach was used to empirically estimate $L_0$, $Q$, and $\Delta$. Details of the neural network architecture, the BNB algorithm, and the sampling approach can be found in Appendix F.\\ 

The constants $L_0$ and $Q$ are non-decreasing in $\kappa$, but it was observed that our estimates can be decreased significantly by decreasing $\kappa$ without having much of an impact on training performance. This enabled reasonable choices for the required number of epochs implied by Corollary \ref{complexity} with a choice of $\kappa=0.2$ and $0.22$ for the MN and FMN datasets. We focused on the first-order version of SPA and ran it according to Corollary \ref{complexity} to achieve an expected $(\epsilon_1,\epsilon_2)$-stationary point for $\epsilon_1=\epsilon_2=1/3$. For the MN and FMN datasets, $\rho$ was chosen as $\rho=2.5$ and $2.75$ to minimize the required number of epochs, searching over a grid of 0.25 increments. We want to highlight that these parameters were chosen solely to ensure an adequate number of epochs, and similar or better solutions are expected for larger $\kappa$, smaller $(\epsilon_1,\epsilon_2)$, with reasonable values of $\rho$, e.g. within the domain of Table \ref{t:1}, but will require longer training times.\\ 

Weights and biases were grouped together by filter and neuron for the convolutional and fully connected layers to generate the partition $\{w^i\}_{i=1}^n$, where $n=236$. The penalty $p_i$ for $i\in[n]$ was set to the dimension of each subset $w^i$, $p_i=d_i$. The parameter $m$ of $C$ was chosen as $m=(1-s)d$ where $s\in(0,1)$ is the chosen sparsity level and $d=44426$. Trying different values of $s$ at $0.05$ increments, layer collapse occurred with $w_1=\Pi_{C\cap\overline{B}_{\beta}}(w_0)$ for randomly initialized $w_0$ and $s=0.7$, so we restricted these experiments to $s\leq 0.65$. The values of $\alpha$, $K$, and $M$ were set according to Corollary \ref{complexity}, and $\beta$ was set to $\beta=0.99(\kappa-\alpha/2)$, such that $\kappa>\beta+\frac{\alpha}{2}$ following Section \ref{smoothing}. The value of $\eta$ was set according to Theorem \ref{maincon}. Table \ref{t:2} presents the values of the aforementioned estimated or computed parameters.\\ 

%\begin{table}[h]
\begin{table}
	\caption{Parameters for MN and FMN datasets for $s\in\{0.65,0.5\}$. All parameters were estimated or given explicitly by the theory of the paper given the choices for $\kappa$, $(\epsilon_1,\epsilon_2)$, and $\rho$, except $\beta=0.99(\kappa-\alpha/2)$, which satisfies $\kappa>\beta+\frac{\alpha}{2}$.}
	\label{t:2}
	\vspace{-0.25cm}
	\begin{center}
		\tabulinesep=1mm
		\begin{tabu}{c|cccccccc}
			\hline
			&$L_0$&$Q$&$\Delta$&$\alpha$&$\beta$&$\eta$&$K$&$M$\\
			\hline	
			MN &8.53E-2&7.49E-3&2.31&6.42E-2&1.66E-1&4.76E-4&4.38E5&2\\
			FMN&1.09E-1&1.22E-2&2.30&5.05E-2&1.93E-1&2.68E-4&7.07E5&3\\
			\hline
		\end{tabu}
	\end{center}
\end{table}

%\begin{figure}[H]
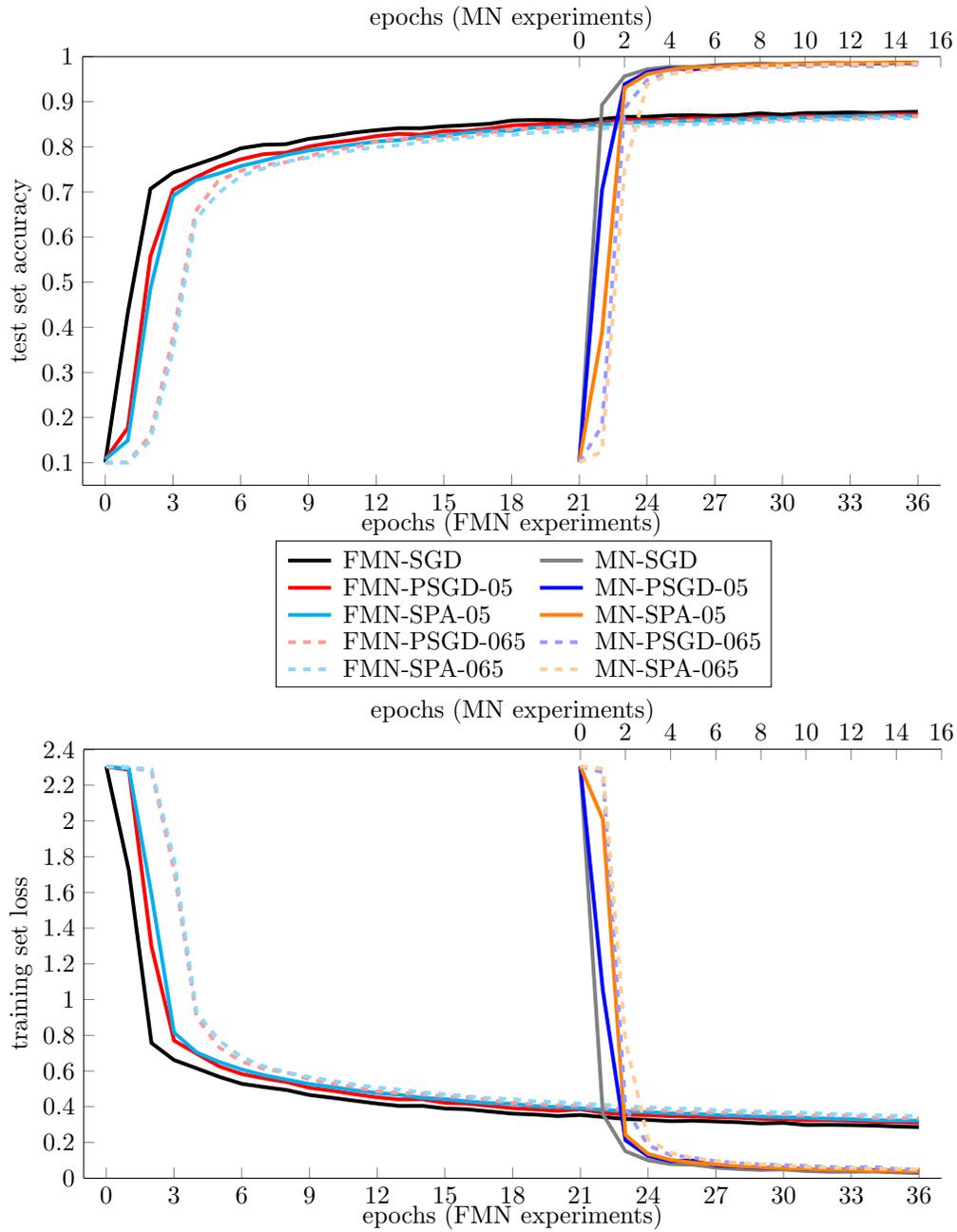
\begin{figure}
	\centering
	\pgfplotsset{width=12cm,height=6cm,compat=1.3}
	\begin{tikzpicture}
		\footnotesize
		\pgfplotsset{scale only axis}
		\begin{axis}[
			axis y line*=left,axis x line*=bottom,xmin=-1,xmax=37,ymin=0.05,ymax=1,
			xtick={0,3,...,36},
			ytick={0,0.1,...,1},
			xlabel=epochs (FMN experiments),
			ylabel=test set accuracy,
			y label style={at={(axis description cs:-0.05,.5)},anchor=south},
			x label style={at={(axis description cs:.5,-0.14)},anchor=south}]
			\addplot[line width=1.5pt,draw=black]
			table[x=x,y=y-test]
			{SGD-FM.dat};\label{plot_three}		
			\addplot[line width=1.5pt,draw=red]
			table[x=x,y=y-test]
			{SGD-FM-05.dat};\label{plot_thirteen}
			\addplot[line width=1.5pt,draw=cyan]
			table[x=x,y=y-test]
			{SPA-FM-05.dat};\label{plot_seven}
			\addplot[line width=1.5pt,draw=red!40!white,dashed]
			table[x=x,y=y-test]
			{SGD-FM-035.dat};\label{plot_eleven}
			\addplot[line width=1.5pt,draw=cyan!40!white,dashed]
			table[x=x,y=y-test]
			{SPA-FM-035.dat};\label{plot_five}			
		\end{axis}
		\begin{axis}[
			axis y line=none,axis x line*=top,ymin=0.05,ymax=1,
			xmin=-22,xmax=16,
			xtick={0,2,...,16},
			xlabel=epochs (MN experiments),
			x label style={at={(axis description cs:.5,1.04)},anchor=south}]			
			\addplot[line width=1.5pt,draw=black!50!white]
			table[x=x,y=y-test]
			{SGD-M.dat};\label{plot_one}
			\addplot[line width=1.5pt,draw=blue]
			table[x=x,y=y-test]
			{SGD-M-05.dat};\label{plot_nineteen}
			\addplot[line width=1.5pt,draw=orange]
			table[x=x,y=y-test]
			{SPA-M-05.dat};\label{plot_nine}
			\addplot[line width=1.5pt,draw=blue!40!white,dashed]
			table[x=x,y=y-test]
			{SGD-M-035.dat};\label{plot_seventeen}
			\addplot[line width=1.5pt,draw=orange!40!white,dashed]
			table[x=x,y=y-test]
			{SPA-M-035.dat};\label{plot_twentyone}		
		\end{axis}	
		\matrix[matrix of nodes,anchor=west,xshift=2.7cm,yshift=-1.8cm,inner sep=0.2em,draw,
		column 2/.style={anchor=base west},
		column 4/.style={anchor=base west}]
		{\ref{plot_three}&FMN-SGD&[1pt]\ref{plot_one}&MN-SGD\\
			\ref{plot_thirteen}&FMN-PSGD-05&[1pt]\ref{plot_nineteen}&MN-PSGD-05\\
			\ref{plot_seven}&FMN-SPA-05&[1pt]\ref{plot_nine}&MN-SPA-05\\
			\ref{plot_eleven}&FMN-PSGD-065&[1pt]\ref{plot_seventeen}&MN-PSGD-065\\
			\ref{plot_five}&FMN-SPA-065&[1pt]\ref{plot_twentyone}&MN-SPA-065\\};
	\end{tikzpicture}
	\pgfplotsset{width=12cm,height=6cm,compat=1.3}
	\begin{tikzpicture}
		\footnotesize
		\pgfplotsset{scale only axis}
		\begin{axis}[
			axis y line*=left,axis x line*=bottom,
			xmin=-1,xmax=37,ymin=0,ymax=2.4,xtick={0,3,...,36},
			ytick={0,0.2,...,2.4},xlabel=epochs (FMN experiments),
			ylabel=training set loss,
			y label style={at={(axis description cs:-0.05,.5)},anchor=south},
			x label style={at={(axis description cs:.5,-0.14)},anchor=south}]
			\addplot[line width=1.5pt,draw=black]
			table[x=x,y=y-train]
			{SGD-FM.dat};\label{plot_four}	
			\addplot[line width=1.5pt,draw=red]
			table[x=x,y=y-train]
			{SGD-FM-05.dat};\label{plot_fourteen}
			\addplot[line width=1.5pt,draw=cyan]
			table[x=x,y=y-train]
			{SPA-FM-05.dat};\label{plot_eight}
			\addplot[line width=1.5pt,draw=red!40!white,dashed]
			table[x=x,y=y-train]
			{SGD-FM-035.dat};\label{plot_twelve}
			\addplot[line width=1.5pt,draw=cyan!40!white,dashed]
			table[x=x,y=y-train]
			{SPA-FM-035.dat};\label{plot_six}
		\end{axis}
		\begin{axis}[
			axis y line=none,axis x line*=top,ymin=0,ymax=2.4,
			xmin=-22,xmax=16,
			xtick={0,2,...,16},
			xlabel=epochs (MN experiments),
			x label style={at={(axis description cs:.5,1.04)},anchor=south}]
			\addplot[line width=1.5pt,draw=gray]
			table[x=x,y=y-train]
			{SGD-M.dat};\label{plot_two}
			\addplot[line width=1.5pt,draw=blue]
			table[x=x,y=y-train]
			{SGD-M-05.dat};\label{plot_twenty}
			\addplot[line width=1.5pt,draw=orange]
			table[x=x,y=y-train]
			{SPA-M-05.dat};\label{plot_ten}
			\addplot[line width=1.5pt,draw=blue!40!white,dashed]
			table[x=x,y=y-train]
			{SGD-M-035.dat};\label{plot_eighteen}
			\addplot[line width=1.5pt,draw=orange!40!white,dashed]
			table[x=x,y=y-train]
			{SPA-M-035.dat};\label{plot_twentytwo}
		\end{axis}
	\end{tikzpicture}
	\caption{Test set accuracy and training set loss of PSGD and SPA. The numbers equal the sparsity level $s$.} \label{T3}
\end{figure}

SPA was compared to projected mini-batch SGD (PSGD) using the same parameterization but with no randomized smoothing, i.e. $u^{k,i}=0\text{ }\forall k,i$, and unconstrained mini-batch SGD (SGD), run identically to PSGD but with no projection. All algorithms were run 3 times for $\lceil KM/\M\rceil$ epochs, where $\M=60,000$ is the training set size, with the output averaged together. The experiments were run in Python 3.6.13 with Pytorch 1.8.1 on a server running Ubuntu 18.04.5 LTS with an Intel Xeon E5-2698 v4 CPU and an Nvidia Titan V GPU. Figure \ref{T3} plots the test set accuracy and the training set loss. The performance of SPA and PSGD are very similar, with better performance in earlier epochs with $s=0.5$, and with all algorithms converging closely to SGD in later epochs.\\ 

The 0-1 knapsack problem did not pose a significant computational bottleneck to the training. An experiment measuring the computation time of $\Pi_{C\cap \overline{B}_{\beta}}(\cdot)$ was conducted for the first 100 projections of 100 trials of the experimental setup of FMN-SGD-065. The projection was found to be the most challenging for $w_0$, with an average computation time of 0.168 seconds, with the remaining projections having an average computation time of 0.0790 seconds. 

\section{Conclusion}
\label{conclusion}
This paper studied theoretical aspects of structured sparsity for deep neural network training. A weighted group $l_0$-norm constraint was proposed and the projection operator and normal cone of this set were presented. The computational complexities of a zeroth and first-order stochastic projection algorithm for constrained Lipschitz continuous loss functions were given for $(\epsilon_1,\epsilon_2)$ and $(\epsilon_3,\epsilon_4)$-stationary points in expectation and high probability, as well as a method with an asymptotic convergence guarantee to a C-M stationary point almost surely.

\appendix
\section*{Appendix A. Section \ref{constraint} Proof}

\begin{proof}{\it (Proposition \ref{closedset})} For a point $x\in \{w^i\neq 0\}$, choosing a $j$ such that $x^i_j\neq 0$, it follows that $B(x,|x^i_j|/2)\in \{w^i\neq 0\}$ as well, proving that $\{w^i\neq 0\}$ is an open set. 
	The lower level sets of $p_i\ind_{\{w^i\neq 0\}}(w)$,
	\begin{alignat}{6}
		&\{w\in\RR^d:p_i\ind_{\{w^i\neq 0\}}(w)\leq \lambda\}=
		\begin{cases}
			\RR^d\quad\text{for } p_i\leq \lambda\\
			\{w^i\neq 0\}^c\quad\text{for } 0\leq\lambda<p_i\\
			\{\emptyset\}\quad\text{for } \lambda<0\\
		\end{cases}\nonumber
	\end{alignat}
	are closed for all $\lambda\in\RR$, which holds if and only if $p_i\ind_{\{w^i\neq 0\}}(w)$ is a lower semicontinuous function \citep[Theorem 1.6]{rockafellar2009}. A function $h$ is lower semicontinuous at $\bar{w}\in\RR^d$ if $\liminf_{w\rightarrow \bar{w}}h(w)\geq h(\bar{w})$ \citep[Definition 1.5]{rockafellar2009}. For all $\bar{w}\in\RR^d$,
	\begin{alignat}{6}
		\liminf_{w\rightarrow \bar{w}} \sum_{i=1}^np_i\ind_{\{w^i\neq 0\}}(w)&\geq \sum_{i=1}^n\liminf_{w\rightarrow \bar{w}} p_i\ind_{\{w^i\neq 0\}}(w)\nonumber\\
		&\geq \sum_{i=1}^np_i\ind_{\{\bar{w}^i\neq 0\}}(\bar{w}),\nonumber 
	\end{alignat} 
	where the second inequality uses the lower semicontinuity of each $p_i\ind_{\{\bar{w}^i\neq 0\}}(\bar{w})$. The lower semicontinuity of $\sum_{i=1}^np_i\ind_{\{w^i\neq 0\}}(w)$ proves that its lower level set $C$ is closed.
\end{proof}

\section*{Appendix B. Section \ref{smoothing} Proofs}

\begin{proof}{\it (Proposition \ref{lipschitz_f})} Given that $\EE[L_0(\xi)^2]<\infty$ it holds that $L_0(\xi)$ is integrable, $\EE[L_0(\xi)]<\infty$ \citep[Proposition 6.12]{folland1999}. For any $w,w'\in \overline{B}_{\kappa}$,
	\begin{alignat}{6}
		|f(w)-f(w')|&&=&|\EE[F(w,\xi)-F(w',\xi)]|\nonumber\\
		&&\leq&\EE[|F(w,\xi)-F(w',\xi)|]\nonumber\\
		&&\leq&\EE[L_0(\xi)||w-w'||_2]\nonumber\\
		&&=&L_0||w-w'||_2,\nonumber
	\end{alignat}
	where the first inequality uses Jensen's inequality and the second inequality holds given that \eqref{loclip} holds for almost all $\xi\in\RR^p$.
\end{proof}

\begin{proof}{\it (Proposition \ref{grad})} Let $\overline{w}\in \overline{B}_{\beta}$, let $\{h_j\}\subset \RR$ be a sequence such that $0<|h_j|\leq \kappa-(\beta+\frac{\alpha}{2})$ with  $\lim\limits_{j\rightarrow\infty}h_j\rightarrow 0$, and let $e_i$ equal the $i^{th}$ unit coordinate vector. Given the Lipschitz continuity of $f(w)$ over $\overline{B}_{\kappa}$, 
	\begin{alignat}{6}
		\lim_{j\rightarrow \infty} h_j^{-1}(f(\overline{w}+u+h_je_i)-f(\overline{w}+u))=\widetilde{\nabla}_i f(\overline{w}+u)\nonumber
	\end{alignat}	
	for almost every $u\in \overline{B}_{\alpha/2}$. For all $j\in \NN$, 
	\begin{alignat}{6}
		|h_j^{-1}||f(\overline{w}+u+h_je_i)-f(\overline{w}+u)|\leq L_0|h_j^{-1}||h_je_i|=L_0\in L^1(P_u).\nonumber
	\end{alignat}
	Applying the dominated convergence theorem,
	\begin{alignat}{6}
		\EE\bigg[\widetilde{\nabla}_i f(\overline{w}+u)\bigg]&=\lim_{j\rightarrow \infty} \EE[h_j^{-1}(f(\overline{w}+u+h_je_i)-f(\overline{w}+u))]\nonumber\\
		&=\lim_{j\rightarrow \infty} h_j^{-1}(f_{\alpha}(\overline{w}+h_je_i)-f_{\alpha}(\overline{w})).\nonumber
	\end{alignat}
	Given that $\EE[|\widetilde{\nabla}_i f(\overline{w}+u)|]\leq L_0$ \citep[Lemma 6]{metel2021}, we can use Fubini's theorem,  
	\begin{alignat}{6}
		\EE\bigg[\widetilde{\nabla}_i f(\overline{w}+u)\bigg]&=\frac{1}{\alpha^d}\int_{-\alpha/2}^{\alpha/2}\cdots\int_{-\alpha/2}^{\alpha/2}\int_{-\alpha/2}^{\alpha/2}\cdots\int_{-\alpha/2}^{\alpha/2}\int_{-\alpha/2}^{\alpha/2}\widetilde{\nabla}_i f(\overline{w}+u)du_idu_1\cdots du_{i-1}du_{i+1}\cdots du_d\nonumber\\
		&=\frac{1}{\alpha^d}\int_{-\alpha/2}^{\alpha/2}\cdots\int_{-\alpha/2}^{\alpha/2}\int_{-\alpha/2}^{\alpha/2}\cdots\int_{-\alpha/2}^{\alpha/2}df_i(\overline{w},u_{\setminus i})du_1\cdots du_{i-1}du_{i+1}\cdots du_d\nonumber\\
		&=\alpha^{-1}\EE[df_i(\overline{w},u_{\setminus i})],\nonumber
	\end{alignat}
	where the second equality uses the fundamental theorem of calculus for Lebesgue integrals given that ${\widetilde{\nabla}_i f(\overline{w}+u)=\nabla_i f(\overline{w}+u)}$ for almost all $u\in [-\alpha/2,\alpha/2]^d$: For a fixed 
	$u'_{\setminus i}\in [-\alpha/2,\alpha/2]^{d-1}$, let 
	$u'(u_i):=[u'_1,...,u'_{i-1},u_i,u'_{i+1},...,u'_d]^T$ for $u_i\in [-\alpha/2,\alpha/2]$.
	Assume that $u'_{\setminus i}$ is chosen such that $f(\overline{w}+u'(u_i))$ is differentiable with  
	$\widetilde{\nabla}_i f(\overline{w}+u'(u_i))=\nabla_i f(\overline{w}+u'(u_i))$ for almost all $u_i\in [-\alpha/2,\alpha/2]$. Given that $f(\overline{w}+u'(u_i))$ is absolutely continuous for $u_i\in [-\alpha/2,\alpha/2]$, ${\int_{-\alpha/2}^{\alpha/2}\widetilde{\nabla}_i f(\overline{w}+u'(u_i))du_i=df_i(\overline{w},u'_{\setminus i})}$. Since this holds for almost all $u'_{\setminus i}\in [-\alpha/2,\alpha/2]^{d-1}$, the second equality holds.\\
	
	The equality $\EE[\widetilde{\nabla} f(\overline{w}+u)]=\EE[\widetilde{\nabla} F(\overline{w}+u,\xi)]$ holds by \citep[Property 2]{metel2021}, where it is shown that for almost all $u\in\overline{B}_{\alpha/2}$, $\widetilde{\nabla} f(\overline{w}+u)=\EE_{\xi}[\widetilde{\nabla} F(\overline{w}+u,\xi)]$ in our problem setting.\\ 
	
	As $\overline{w}\in \overline{B}_{\beta}$ and the sequence $\{h_i\}\subset \overline{B}_{\kappa-(\beta+\frac{\alpha}{2})}$ were arbitrary, for all $w\in \overline{B}_{\beta}$,
	\begin{alignat}{6}
		\alpha^{-1}\EE[df_i(w,u_{\setminus i})]&=\EE[\widetilde{\nabla} F_i(w+u,\xi)]\nonumber\\
		&=\EE[\widetilde{\nabla}_i f(w+u)]\nonumber\\
		&=\lim_{h\rightarrow 0} h^{-1}(f_{\alpha}(w+he_i)-f_{\alpha}(w))\nonumber\\
		&=\frac{\partial f_{\alpha}}{\partial w_i}(w)\nonumber\\
		&=\nabla_i f_{\alpha}(w),\label{gradd}
	\end{alignat}
	where the third equality holds using the sequential criterion of a limit, and the last equality holds given that the partial derivatives are (Lipschitz) continuous: For any $w,w'\in \overline{B}_{\beta}$,
	\begin{alignat}{6}
		&|\alpha^{-1}\EE[df_i(w,u_{\setminus i})]-\alpha^{-1}\EE[df_i(w',u_{\setminus i})]|\nonumber\\
		\leq & \alpha^{-1}\EE[|df_i(w,u_{\setminus i})-df_i(w',u_{\setminus i})|]\nonumber\\
		\leq & 2\alpha^{-1}L_0||w-w'||_2,\nonumber
	\end{alignat}
	where a proof of the last inequality can be found at \eqref{bound}. Given that for any $w\in \overline{B}_{\beta}$,
	\begin{alignat}{6}
		\alpha^{-1}\EE[|dF_i(w,u_{\setminus i},\xi)|]
		\leq\alpha^{-1}\EE[L_0(\xi)\alpha]
		=L_0,\nonumber
	\end{alignat}	
	Fubini's theorem can be applied: 
	\begin{alignat}{6}
		&\alpha^{-1}\EE[dF_i(w,u_{\setminus i},\xi)]
		&&=\alpha^{-1}\EE_{u_{\setminus i}}[\EE_\xi[dF_i(w,u_{\setminus i},\xi)]]\nonumber\\
		&&&=\alpha^{-1}\EE[df_i(w,u_{\setminus i})]\nonumber\\
		&&&=\nabla_i f_{\alpha}(w)\nonumber
	\end{alignat}
	from \eqref{gradd}.
\end{proof}

\begin{proof}{\it (Proposition \ref{clarkeps})} Given that $u\in \overline{B}_{\frac{\alpha}{2}}$, $w+u\in \overline{B}(w,\sqrt{d}\frac{\alpha}{2})\subset \overline{B}_{\kappa}$, hence 
	$\overline{\partial}_{\widehat{\alpha}} f(w)=\text{co}\{\overline{\partial} f(x): x\in \overline{B}(w,\widehat{\alpha})\}$ is well defined. The gradient 
	$\nabla f(w+u)\in \overline{\partial} f(w+u)$ wherever $f(w+u)$ is differentiable \citep[Proposition 2.2.2]{clarke1990}, which is for almost all $u\in \overline{B}_{\frac{\alpha}{2}}$. It follows that for almost all $u\in \overline{B}_{\frac{\alpha}{2}}$, $\widetilde{\nabla} f(w+u)\in  \overline{\partial}_{\widehat{\alpha}}f(w)$ and 
	$\nabla f_{\alpha}(w)=\EE[\widetilde{\nabla} f(w+u)]\in  \overline{\partial}_{\widehat{\alpha}}f(w)$.
\end{proof}

\begin{proof}{\it (Proposition \ref{lipcont})}
	\text{ }
	\begin{enumerate}
		\item For all $w,w'\in \overline{B}_{\beta}$,
		\begin{alignat}{6}
			|f_{\alpha}(w)-f_{\alpha}(w')|&=|\EE[f(w+u)-f(w'+u)]|\leq L_0||w-w'||_2.\nonumber
		\end{alignat}
		\item For all $w\in \overline{B}_{\beta}$,
		\begin{alignat}{6}
			|f_{\alpha}(w)-f(w)|&=|\EE[f(w+u)-f(w)]|\leq L_0\EE[||u||_2]\leq L_0\sqrt{\frac{d}{12}}a,\nonumber
		\end{alignat}	
		where the expected distance from $u$ to the origin is bounded by  
		\begin{alignat}{6}
			\EE[\norm{u}_2]&\leq \sqrt{\EE[\sum_{i=1}^du_i^2]}&= \sqrt{d}\sqrt{\EE[u_i^2]}&=\sqrt{d}\frac{\alpha}{\sqrt{12}}.\nonumber
		\end{alignat}
		\item This can be proven by contradiction.
		Assuming that $f_{\alpha}(w_{\alpha}^*)-f(w^*)>\alpha L_0\sqrt{\frac{d}{12}}$ and given that ${f_{\alpha}(w^*)-f(w^*)\leq\alpha L_0\sqrt{\frac{d}{12}}}$ from statement 2, 
		\begin{alignat}{6}
			&&f_{\alpha}(w_{\alpha}^*)&>\alpha L_0\sqrt{\frac{d}{12}}+f(w^*)\geq f_{\alpha}(w^*),\nonumber
		\end{alignat}
		contradicting the optimality of $w_{\alpha}^*$ for $f_{\alpha}(w)$. 
		Similarly if $f(w^*)-f_{\alpha}(w_{\alpha}^*)>\alpha L_0\sqrt{\frac{d}{12}}$, statement 2 
		gives ${f(w_{\alpha}^*)-f_{\alpha}(w_{\alpha}^*)\leq \alpha L_0\sqrt{\frac{d}{12}}}$, 
		from which
		\begin{alignat}{6}
			&&f(w^*)&>\alpha L_0\sqrt{\frac{d}{12}}+f_{\alpha}(w_{\alpha}^*)\geq f_{\alpha}(w^*),\nonumber
		\end{alignat}
		contradicting the optimality of $w^*$ for $f(w)$.
		\item For any two values $w,w'\in \overline{B}_{\beta}$,  	
		\begin{alignat}{6}
			|f_{\alpha}(w)-f_{\alpha}(w')|&=|\EE[f(w+u)]-\EE[f(w'+u)]|\nonumber\\
			&=|\EE[f(w+u)-f(w'+u)]|\nonumber\\
			&\leq\EE[L_0||w-w'||_2]\nonumber\\
			&\leq L_0\sqrt{\sum_{i=1}^d(2\beta)^2}\nonumber\\
			&=L_0\sqrt{d}2\beta.\nonumber
		\end{alignat}
	\end{enumerate}
\end{proof}

\begin{proof}{\it (Proposition \ref{gradcont})} From Proposition \ref{grad}, 
	\begin{alignat}{6}
		||\nabla f_{\alpha}(w)-\nabla f_{\alpha}(w')||_2=&\alpha^{-1}||\EE[df(w,u)-df(w',u)]||_2\nonumber\\
		=&\alpha^{-1}\sqrt{\sum_{i=1}^d(\EE[df_i(w,u_{\setminus i})-df_i(w',u_{\setminus i})])^2}.\label{norm}
	\end{alignat}
	Focusing on $df_i(w,u_{\setminus i})-df_i(w',u_{\setminus i})$,
	\begin{alignat}{6}
		&|df_i(w,u_{\setminus i})-df_i(w',u_{\setminus i})|\nonumber\\
		=&|f(w_1+u_1,...,w_{i-1}+u_{i-1},w_i+\frac{\alpha}{2},w_{i+1}+u_{i+1},...,w_d+u_d)\nonumber\\
		-&f(w_1+u_1,...,w_{i-1}+u_{i-1},w_i-\frac{\alpha}{2},w_{i+1}+u_{i+1},...,w_{d}+u_d)\nonumber\\
		-&f(w'_1+u_1,...,w'_{i-1}+u_{i-1},w'_i+\frac{\alpha}{2},w'_{i+1}+u_{i+1},...,w'_d+u_d)\nonumber\\
		+&f(w'_1+u_1,...,w'_{i-1}+u_{i-1},w'_i-\frac{\alpha}{2},w'_{i+1}+u_{i+1},...,w'_{d}+u_d)|\nonumber\\
		\leq&|f(w_1+u_1,...,w_{i-1}+u_{i-1},w_i+\frac{\alpha}{2},w_{i+1}+u_{i+1},...,w_d+u_d)\nonumber\\
		-&f(w'_1+u_1,...,w'_{i-1}+u_{i-1},w'_i+\frac{\alpha}{2},w'_{i+1}+u_{i+1},...,w'_d+u_d)|\nonumber\\
		+&|f(w'_1+u_1,...,w'_{i-1}+u_{i-1},w'_i-\frac{\alpha}{2},w'_{i+1}+u_{i+1},...,w'_{d}+u_d)\nonumber\\
		-&f(w_1+u_1,...,w_{i-1}+u_{i-1},w_i-\frac{\alpha}{2},w_{i+1}+u_{i+1},...,w_{d}+u_d)|\nonumber\\
		\leq& 2L_0||w-w'||_2.\label{bound}
	\end{alignat}
	Plugging \eqref{bound} into \eqref{norm},
	\begin{alignat}{6}
		||\nabla f_{\alpha}(w)-\nabla f_{\alpha}(w')||_2\leq&\alpha^{-1}\sqrt{\sum_{i=1}^d(2L_0||w-w'||_2)^2}\nonumber\\
		=&\alpha^{-1}\sqrt{d}2L_0||w-w'||_2.\nonumber
	\end{alignat}
\end{proof}

\begin{proof}{\it (Proposition \ref{dq})} To streamline the proof, let $G:=\frac{1}{M}\sum_{i=1}^M g^i$, where
	\begin{alignat}{6}
		g^i:=\begin{cases}
			\alpha^{-1}dF(w,u^{i},\xi^{i}) &\text{ for } 1 \text{ and } 3\\ 
			\widetilde{\nabla} F(w+u^{i},\xi^{i})& \text{ for } 2 \text{ and } 4.\\ 
		\end{cases}\nonumber
	\end{alignat}
	Proof of 1 and 2:
	\begin{alignat}{6}
		&&\EE[||\nabla f_{\alpha}(w)-G||^2_2]=&\EE[||\EE[G]-G||^2_2]\nonumber\\
		&&=&\EE\bigg[\sum_{j=1}^d(\EE[G_j]-G_j)^2\bigg]\nonumber\\
		&&=&\EE\bigg[\sum_{j=1}^d(\frac{1}{M}\sum_{i=1}^M(\EE[G_j]-g^i_j))^2\bigg]\nonumber\\
		&&=&\frac{1}{M^2}\EE\bigg[\sum_{j=1}^d\bigg(\sum_{i=1}^M(\EE[G_j]-g^i_j)^2+2\sum_{i=1}^M\sum_{l=1}^{i-1}(\EE[G_j]-g^i_j)(\EE[G_j]-g^l_j)\bigg)\bigg]\nonumber\\
		&&=&\frac{1}{M^2}\sum_{j=1}^d\sum_{i=1}^M\EE[(\EE[G_j]-g^i_j)^2]\nonumber\\
		&&=&\frac{1}{M}\sum_{j=1}^d\EE[(\EE[G_j]-g^i_j)^2]\nonumber\\
		&&=&\frac{1}{M}\sum_{j=1}^d(\EE[G_j]^2-2\EE[G_j]\EE[g^i_j]+\EE[(g^i_j)^2])\nonumber\\
		&&=&\frac{1}{M}\sum_{j=1}^d(-\EE[G_j]^2+\EE[(g^i_j)^2])\nonumber\\
		&&\leq&\frac{1}{M}\sum_{j=1}^d\EE[(g^i_j)^2].\nonumber
	\end{alignat}
	For 1:
	\begin{alignat}{6}
		&&\frac{1}{M}\sum_{j=1}^d\EE[(g^i_j)^2]=&\frac{1}{M}\sum_{j=1}^d\EE[(\alpha^{-1}dF_j(w,u^i
		_{\setminus j},\xi^{i}))^2]\nonumber\\
		&&=&\frac{1}{\alpha^2M}\sum_{j=1}^d\EE[dF_j(w,u^i_{\setminus j},\xi^i)^2]\nonumber\\
		&&\leq&\frac{1}{\alpha^2M}\sum_{j=1}^d\EE[L_0(\xi)^2\alpha^2]\nonumber\\
		&&=&\frac{d}{M}\EE[L_0(\xi)^2].\label{13bound}
	\end{alignat}
	
	For 2:
	\begin{alignat}{6}
		&&\frac{1}{M}\sum_{j=1}^d\EE[(g^i_j)^2]=&\frac{1}{M}\sum_{j=1}^d\EE[\widetilde{\nabla} F(w+u^{i},\xi^{i})^2]\nonumber\\
		&&=&\frac{1}{M}\EE[||\widetilde{\nabla} F(w+u^{i},\xi^{i})||^2_2]\nonumber\\
		&&\leq&\frac{1}{M}\EE[L_0(\xi)^2],\label{24bound}
	\end{alignat}
	where the last inequality follows from \citep[Lemma 6]{metel2021}.\\
	
	Proof of 3 and 4:
	\begin{alignat}{6}
		&&\EE[||G||^2_2]=&\EE[||\frac{1}{M}\sum_{i=1}^Mg^i||^2_2]\nonumber\\
		&&=&\EE[\sum_{j=1}^d(\frac{1}{M}\sum_{i=1}^Mg^i_j)^2]\nonumber\\
		&&\leq&\EE[\sum_{j=1}^d\frac{1}{M}\sum_{i=1}^M(g^i_j)^2]\nonumber\\
		&&=&\EE[\sum_{j=1}^d(g^i_j)^2],\nonumber
	\end{alignat}
	where the inequality uses Jensen's inequality. For 3, from \eqref{13bound}:
	\begin{alignat}{6}
		&&\EE[\sum_{j=1}^d(g^i_j)^2]\leq&d\EE[L_0(\xi^{i})^2],\nonumber
	\end{alignat}
	and for 4, from \eqref{24bound}:
	\begin{alignat}{6}
		&&\EE[\sum_{j=1}^d(g^i_j)^2]\leq&\EE[L_0(\xi^{i})^2].\nonumber
	\end{alignat}
\end{proof}

\section*{Appendix C. Section \ref{constraintprop} Proofs}

\begin{proof}{\it (Proposition \ref{projection})} Following \citep[Chapter 1.G.]{rockafellar2009}, the projection of a point $w$ onto $C\cap \overline{B}_{\beta}$ can be written as 
	\begin{alignat}{6}
		\Pi_{C\cap \overline{B}_{\beta}}(w)=\argmin_{x\in\RR^d}\text{ }&\{\delta_{C\cap \overline{B}_{\beta}}(x)+\sum_{i=1}^n||x^i-w^i||^2_2\}.\nonumber
	\end{alignat}
	For the squared distance function of a point $w\in\RR^d$ from $C\cap \overline{B}_{\beta}$, 
	\begin{alignat}{6}
		d^2_{C\cap \overline{B}_{\beta}}(w)=\inf_{x\in\RR^d}\text{ }&\{\delta_{C\cap \overline{B}_{\beta}}(x)+\sum_{i=1}^n||x^i-w^i||^2_2\},\nonumber
	\end{alignat}
	let $x^*$ be a minimizer. If $(x^*)^i\neq 0$ then it will be the optimal solution of 
	\begin{alignat}{6}
		&&\min\limits_{x\in\RR^{d_i}}\text{ }&||x^i-w^i||^2_2\nonumber\\
		&&\st\text{ }&|x_j^i|\leq \beta\quad j=1,2,...,d_i,\nonumber
	\end{alignat}
	equal to $(x^*)_j^i=\sgn(w_j^i)\min(|w^i_j|,\beta)$ for $j=1,2,...,d_i$, with the optimal objective value of $||\max(|w^i|-\beta,0)||^2_2$. The function $d^2_{C\cap \overline{B}_{\beta}}$ can then be written as
	\begin{alignat}{6}
		d^2_{C\cap \overline{B}_{\beta}}(w)=\min_{z\in\{0,1\}^n}\text{ }&\sum_{i=1}^nz_i||\max(|w^i|-\beta,0)||^2_2+(1-z_i)||w^i||^2_2\nonumber\\
		\text{s.t. }&\sum_{i=1}^nz_ip_i\leq m,\nonumber
	\end{alignat}	
	where each $z_i$ decides if $x^i$ is allowed to be non-zero. This optimization problem has the same set of optimal solutions as  
	\begin{alignat}{6}	
		\min_{z\in\{0,1\}^n}\text{ }&\sum_{i=1}^nz_i(||\max(|w^i|-\beta,0)||^2_2-||w^i||^2_2)\nonumber\\
		\text{s.t. }&\sum_{i=1}^nz_ip_i\leq m,\nonumber
	\end{alignat}
	and
	\begin{alignat}{6}			
		\max_{z\in\{0,1\}^n}\text{ }&\sum_{i=1}^nz_i(||w^i||^2_2-||\max(|w^i|-\beta,0)||^2_2)\nonumber\\
		\text{s.t. }&\sum_{i=1}^nz_ip_i\leq m,\nonumber
	\end{alignat}
	written in the form of a maximization to match the standard format of knapsack problems.
\end{proof}

\begin{proof}{\it (Proposition \ref{frechnorm})}
	Let $RHS$ equal the right-hand side of \eqref{freeq} and for simplicity let $I:= I(w)$ and $J:= J(w)$. If $w=0$, then $RHS=0\in\widehat{N}(w,C\cap \overline{B}_{\beta})$. Assuming $w\neq 0$, let $v\in RHS$ and $\hat{m}:=\min_{i\in I} ||w^i||_2/2$. For all $x\in B(w,\hat{m})\cap C\cap \overline{B}_{\beta}=:D$, $I\subseteq I(x)$. Given that $x\in C$ and $I\subseteq I(x)$, $I(x)\subseteq I\cup J$: If not, then there exists a $j\in I(x)\setminus I$ such that $j\notin J$, but this implies that $\sum_{i\in I(x)}p_i +p_j\geq\sum_{i\in I}p_i +p_j>m$, which contradicts that $x\in C$. Given that $I(x)\subseteq I\cup J$ for all $x\in D$, 
	\begin{alignat}{6}\label{ineq}
		\langle v,x-w\rangle=\sum_{i\notin I\cup J}\langle v^i,x^i-w^i\rangle+\sum_{i\in I\cup J}\sum_{j\in [d_i]} v^i_j(x^i_j-w^i_j)\leq 0\end{alignat} 
	since for $i\notin I\cup J$, $x^i=w^i=0$ and for $i\in I\cup J$, given that $x\in \overline{B}_{\beta}$, $v^i_j(x^i_j-w^i_j)\leq 0$ from \eqref{freeq}. From \eqref{ineq},
	\begin{alignat}{6}
		\inf_{\gamma>0}\bigg(\sup_{\substack{{0<||x-w||_2<\gamma}\\x\in C\cap \overline{B}_{\beta}}}\frac{\langle v,x-w\rangle}{||x-w||_2}\bigg)\leq \sup_{\substack{0<||x-w||_2<\hat{m}\\x\in C\cap \overline{B}_{\beta}}}\frac{\langle v,x-w\rangle}{||x-w||_2}\leq 0\nonumber
	\end{alignat}
	and $v\in\widehat{N}(w,C\cap \overline{B}_{\beta})$.\\
	
	For a $v\in \widehat{N}(w,C\cap \overline{B}_{\beta})$, assume there exists an $i\in I\cup J$ and $j\in [d_i]$ for which it does not hold that 
	\begin{alignat}{6}
		v^i_j\in
		\begin{cases}
			\RR_{\geq 0}&\text{if }w^i_j=\beta\\
			0&\text{if }|w^i_j|<\beta\\
			\RR_{\leq 0}&\text{if }w^i_j=-\beta.
		\end{cases}\nonumber		
	\end{alignat} 
	
	Consider the sequence $\{x_k\}$ with elements equal to   
	\begin{alignat}{6}
		(x^l_m)_k=\begin{cases}
			w^l_m+\epsilon_kv^l_m & \text{ if } l=i \text{ and } m=j\\
			w^l_m& \text{ otherwise }\\
		\end{cases}\nonumber
	\end{alignat} 
	with 
	\begin{alignat}{6}
		\begin{cases}
			\vspace{2pt}
			0<\epsilon_k\leq -2\frac{\beta}{v^i_j}&\text{ if } v^i_j<0 \text{ and } w^i_j=\beta\\
			0<\epsilon_k\leq \frac{\beta-|w^i_j|}{|v^i_j|}&\text{ if } v^i_j\neq 0 \text{ and } |w^i_j|<\beta\\
			0<\epsilon_k\leq 2\frac{\beta}{v^i_j}&\text{ if } v^i_j>0 \text{ and } w^i_j=-\beta,\\
		\end{cases}\nonumber
	\end{alignat} 
	to ensure that $\{x_k\}\subset\overline{B}_{\beta}$, and assume that $\epsilon_k\rightarrow 0$. It holds that $x_k\xrightarrow{C\cap \overline{B}_{\beta}}w$, since $I(x_k)\subseteq I\cup i$ for all $k\in \NN$. For any $\gamma>0$, there exists a $K_{\gamma}\in \NN$ such that for $k>K_{\gamma}$, $0<||x_k-w||_2<\gamma$, hence
	\begin{alignat}{6}
		\sup_{\substack{{0<||x-w||_2<\gamma}\\x\in C\cap \overline{B}_{\beta}}}\frac{\langle v,x-w\rangle}{||x-w||_2}\geq \frac{\langle v,x_k-w\rangle}{||x_k-w||_2}=\frac{\epsilon_k(v^i_j)^2}{\epsilon_k|v^i_j|}=|v^i_j|>0,\nonumber
	\end{alignat}
	which contradicts that $v\in \widehat{N}(w,C\cap \overline{B}_{\beta})$.   
\end{proof}

\begin{proof}{\it (Proposition \ref{limitnorm})}
	Let $RHS$ equal the right-hand side of \eqref{moreq} and let $I:= I(w)$ and $J:= J(w)$. For any $v\in N(w,C\cap \overline{B}_{\beta})$, there exists sequences $x_k\xrightarrow{C\cap \overline{B}_{\beta}}w$ and $v_k\rightarrow v$ with $v_k\in \widehat N(x_k,C\cap \overline{B}_{\beta})$ for all $k\in \NN$. We first want to show that for $k\in \NN$ sufficiently large, there exist $X(v_k)\in Y$ such that $X(v_k)\subseteq I(x_k)\cup J(x_k)$, implying that $v_k\in RHS$.\\
	
	Assuming $w\neq 0$, there exists an $N\in \NN$ such that for all $k>N$, $\max_{i\in [n]}||x_k^i-w^i||_2<\min_{i\in I} ||w^i||_2/2$, implying that $I\subseteq I(x_k)$. Given that $I\subseteq I(x_k)$, there exists an $X(v_k)\in Y$ such that $X(v_k)\subseteq I(x_k)\cup J(x_k)$: 
	We can choose $X(v_k)$ such that $I(x_k)\subseteq X(v_k)$ given that $I\subseteq I(x_k)$ and $x_k\in C$. If $X(v_k)\not\subset I(x_k)\cup J(x_k)$, then there exists a $j\in X(v_k)\setminus I(x_k)$ such that $j\notin J(x_k)$, but this implies that $\sum_{i\in X(v_k)}p_i\geq\sum_{i\in I(x_k)}p_i +p_j>m$, which contradicts that $X(v_k)\in Y$.\\ 
	
	For the case when $w=0$, $I=\{\emptyset\}$, hence we can choose $X(v_k)$ such that $I(x_k)\subseteq X(v_k)$, from which it must hold again that $X(v_k)\subseteq I(x_k)\cup J(x_k)$. Given the existence of an $X(v_k)\in Y$ such that $X(v_k)\subseteq I(x_k)\cup J(x_k)$, $v_k\in RHS$ for all $k>N$, from which it follows that $v\in RHS$ given that RHS is a closed set:
	For any $w\in \RR^d$, 
	\begin{alignat}{6}
		N^i_j(w):=
		\left\{y\in\RR^d: y^i_j\in\begin{cases}
			\RR_{\geq 0}&\text{if }w^i_j=\beta\\
			0&\text{if }|w^i_j|<\beta\\
			\RR_{\leq 0}&\text{if }w^i_j=-\beta
		\end{cases}\right\}\nonumber
	\end{alignat}
	is a closed set, the intersection of closed sets
	\begin{alignat}{6}\label{capp}
		\bigcap_{\substack{\forall i\in X\\\forall j\in[d_i]}}N^i_j(w)
	\end{alignat}
	for an $X\in Y$ is closed, and the finite union of sets \eqref{capp} over $X\in Y$ 
	\begin{alignat}{6}
		RHS=\bigcup_{X\in Y}\bigcap_{\substack{\forall i\in X\\\forall j\in[d_i]}}N^i_j(w)\nonumber
	\end{alignat}
	is closed.\\
	
	Let $v\in RHS$ and $X(v)\in Y$ such that $\forall i\in X(v)$ and $\forall j\in[d_i]$, 
	\begin{alignat}{6}
		v^i_j\in\begin{cases}
			\RR_{\geq 0}&\text{if }w^i_j=\beta\\
			0&\text{if }|w^i_j|<\beta\\
			\RR_{\leq 0}&\text{if }w^i_j=-\beta.\nonumber
		\end{cases}
	\end{alignat}
	
	Consider the sequence $\{x_k\}$ equal to   
	\begin{alignat}{6}
		x^i_k=\begin{cases}
			w^i& \text{ if } i\in I\\
			\epsilon^i_k & \text{ if } i\in X(v)\setminus I\\
			0 & \text{ otherwise,}\\
		\end{cases}\nonumber
	\end{alignat} 
	where for a $j\in [d_i]$, $0<(\epsilon^i_j)_k\leq \beta$ and $(\epsilon^i_j)_k\rightarrow 0$, and $(\epsilon^i_l)_k=0$ for all $l\in[d_i]\setminus\{j\}$ and $k\in\NN$. It holds that 
	$x_k\xrightarrow{C\cap \overline{B}_{\beta}}w$, the sets $I(x_k)=X(v)$ and $J(x_k)=\{\emptyset\}$, hence $v\in \widehat{N}(x_k,C)$ for all $k\in N$, proving that $v\in N(w,C)$. We note that for the case $X(v)=I$, $v\in \widehat{N}(w,C)\subseteq N(w,C)$ \citep[Eq. (1.6)]{mordukhovich2013} 
\end{proof}

\section*{Appendix D. Section \ref{trainalg} Proofs}

\begin{proof}{\it (Proposition \ref{meas})}
	The measurability of $\Pi_{S\cap \overline{B}_{\beta}}(\cdot)$ follows from \citep[Exercise 14.17 (b)]{rockafellar2009}. In particular, the set $S\cap \overline{B}_{\beta}$ is closed and can be considered as a constant (measurable) set valued function $S\cap \overline{B}_{\beta}:\RR^d\rightrightarrows \RR^d$. Given that the projection operator is a closed (in particular a nonempty and compact) set-valued mapping \citep[Exercise 1.20]{rockafellar2009}, there exists a measurable selection \citep[Corollary 14.6]{rockafellar2009}.\\
	
	The restriction that $\epsilon\leq\frac{\alpha}{2}$ is to ensure that $\overline{\partial}_{\epsilon}f(w)$ is well defined given that it is only assumed that $\kappa>\beta+\frac{\alpha}{2}$.
	The set valued mapping $\overline{\partial}_{\epsilon}f(w)$ is outer semicontinuous for $w\in \overline{B}_{\beta}$, see  \citep[Lemma 2.6]{goldstein1977} for a proof that \citep[Definition 5.4]{rockafellar2009} holds, as is $N(w,S\cap \overline{B}_{\beta})$ for $w\in\RR^d$ given that ${S\cap\overline{{B}}_{\beta}}$ is closed \citep[Theorem 5.7 (a); Page 202]{rockafellar2009}. Since $\overline{\partial}_{\epsilon}f(w)$ is bounded, given that $\overline{\partial}_{\epsilon}f(w)\subseteq\overline{B}(0,L_0)$ \citep[Proposition 2.1.2 (a)]{clarke1990}, $\overline{\partial}_{\epsilon}f(w)+N(w,S\cap \overline{B}_{\beta})$ is outer semicontinuous for $w\in \overline{B}_{\beta}$ \citep[Proposition 5.51 (b)]{rockafellar2009}. The function 
	$\dist(0,\overline{\partial}_{\epsilon}f(w)+N(w,S\cap {\overline{B}}_{\beta}))$ is then lower semicontinuous for $w\in \overline{B}_{\beta}$ \citep[Proposition 5.11]{rockafellar2009}, hence Borel measurable. The Borel measurability of $\dist(0,\nabla f_{\alpha}(w)+N(w,S\cap \overline{B}_{\beta}))$ holds since $f_{\alpha}(w)$ is Lipschitz continuous (Proposition \ref{lipcont}.1), continuously differentiable (Proposition \ref{gradcont}), hence $\nabla f_{\alpha}(w)=\overline{\partial}f_\alpha(w)$ \citep[Page 10]{clarke1990}.  
\end{proof}

\begin{proof}{\it (Theorem \ref{maincon})}
	For simplicity let 
	\begin{alignat}{6}
		G^k:=\begin{cases}
			\frac{1}{M\alpha}\sum_{i=1}^MdF(w^k,u^{k,i},\xi^{k,i}) & \text{ if using } (1)\\ 
			\frac{1}{M}\sum_{i=1}^M\widetilde{\nabla} F(w^k+u^{k,i},\xi^{k,i}) & \text{ if using } (2)\\ 
		\end{cases}\nonumber
	\end{alignat}
	in SPA, where we assume $k\in[K]$. Given that 
	\begin{alignat}{6}
		&&w^{k+1}\in\argmin_{x\in\RR^d}\text{ }\{\delta_{S\cap \overline{B}_{\beta}}(x)+||x-(w^k-\eta G^k)||^2_2\}\label{proopt}
	\end{alignat}
	and $||x-(w^k-\eta G^k)||^2_2$ is differentiable \citep[Theorem 6.12]{rockafellar2009},
	\begin{alignat}{6}
		&&-2\eta G^k-2(w^{k+1}-w^k)&\in&&N(w^{k+1},S\cap \overline{B}_{\beta})\nonumber\\
		\Longrightarrow&&-G^k-\eta^{-1}(w^{k+1}-w^k)&\in&&N(w^{k+1},S\cap \overline{B}_{\beta})\nonumber\\
		\Longrightarrow&&\nabla f_{\alpha}(w^{k+1})-G^k-\eta^{-1}(w^{k+1}-w^k)&\in&&\nabla f_{\alpha}(w^{k+1})+N(w^{k+1},S\cap \overline{B}_{\beta}),\nonumber
	\end{alignat}
	where the second inclusion holds since $N(w^{k+1},S\cap \overline{B}_{\beta})$ is a cone. The final inclusion implies that   
	\begin{alignat}{6}
		&&&\dist(0,\nabla f_{\alpha}(w^{k+1})+N(w^{k+1},S\cap \overline{B}_{\beta}))^2\nonumber\\
		&&\leq& ||\nabla f_{\alpha}(w^{k+1})-G^k-\eta^{-1}(w^{k+1}-w^k)||^2_2\nonumber\\
		&&=& ||\nabla f_{\alpha}(w^{k+1})-G^k||^2_2-2\eta^{-1}\langle \nabla f_{\alpha}(w^{k+1})-G^k,w^{k+1}-w^k\rangle+\eta^{-2}||w^{k+1}-w^k||^2_2\nonumber\\
		&&=& ||\nabla f_{\alpha}(w^{k+1})-\nabla f_{\alpha}(w^k)+\nabla f_{\alpha}(w^k)-G^k||^2_2\nonumber\\
		&&&-2\eta^{-1}\langle \nabla f_{\alpha}(w^{k+1})-G^k,w^{k+1}-w^k\rangle+\eta^{-2}||w^{k+1}-w^k||^2_2\nonumber\\
		&&\leq& \frac{3}{2}||\nabla f_{\alpha}(w^{k+1})-\nabla f_{\alpha}(w^k)||^2_2+3||\nabla f_{\alpha}(w^k)-G^k||^2_2\nonumber\\
		&&&-2\eta^{-1}\langle \nabla f_{\alpha}(w^{k+1})-G^k,w^{k+1}-w^k\rangle+
		\eta^{-2}||w^{k+1}-w^k||^2_2\nonumber\\
		&&\leq& 6\alpha^{-2}dL_0^2||w^{k+1}-w^k||^2_2+3||\nabla f_{\alpha}(w^k)-G^k||^2_2\nonumber\\
		&&&-2\eta^{-1}\langle \nabla f_{\alpha}(w^{k+1})-G^k,w^{k+1}-w^k\rangle+
		\eta^{-2}||w^{k+1}-w^k||^2_2\nonumber\\
		&&=& (6\alpha^{-2}dL_0^2+\eta^{-2})||w^{k+1}-w^k||^2_2+3||\nabla f_{\alpha}(w^k)-G^k||^2_2\nonumber\\
		&&&-2\eta^{-1}\langle \nabla f_{\alpha}(w^{k+1})-\nabla f_{\alpha}(w^k)+\nabla f_{\alpha}(w^k)-G^k,w^{k+1}-w^k\rangle\nonumber\\
		&&\leq& (6\alpha^{-2}dL_0^2+\eta^{-2})||w^{k+1}-w^k||^2_2+3||\nabla f_{\alpha}(w^k)-G^k||^2_2\nonumber\\
		&&&+4\alpha^{-1}\sqrt{d}L_0\eta^{-1}||w^{k+1}-w^k||^2_2-2\eta^{-1}\langle \nabla f_{\alpha}(w^k)-G^k,w^{k+1}-w^k\rangle\nonumber\\
		&&=& (6\alpha^{-2}dL_0^2+4\alpha^{-1}\sqrt{d}L_0\eta^{-1}+\eta^{-2})||w^{k+1}-w^k||^2_2+3||\nabla f_{\alpha}(w^k)-G^k||^2_2\nonumber\\
		&&&-2\eta^{-1}\langle \nabla f_{\alpha}(w^k)-G^k,w^{k+1}-w^k\rangle,\label{join0}
	\end{alignat}
	where the second inequality uses Young's inequality:
	\begin{alignat}{6}
		&&&||\nabla f_{\alpha}(w^{k+1})-\nabla f_{\alpha}(w^k)+\nabla f_{\alpha}(w^k)-G^k||^2_2\nonumber\\
		&&=&||\nabla f_{\alpha}(w^{k+1})-\nabla f_{\alpha}(w^k)||^2_2+
		2\langle \frac{1}{\sqrt{2}}(\nabla f_{\alpha}(w^{k+1})-\nabla f_{\alpha}(w^k)),\sqrt{2}(\nabla f_{\alpha}(w^k)-G^k)\rangle\nonumber\\
		&&+&||\nabla f_{\alpha}(w^k)-G^k||^2_2\nonumber\\
		&&\leq&\frac{3}{2}||\nabla f_{\alpha}(w^{k+1})-\nabla f_{\alpha}(w^k)||^2_2
		+3||\nabla f_{\alpha}(w^k)-G^k||^2_2,\nonumber
	\end{alignat}
	and the third and fourth inequalities use the smoothness of $f_{\alpha}(w)$. By the optimality of $w^{k+1}$ in \eqref{proopt},
	\begin{alignat}{6}
		&&\delta_{S\cap \overline{B}_{\beta}}(w^{k+1})+||w^{k+1}-(w^k-\eta G^k)||^2_2&\leq&&\delta_{S\cap \overline{B}_{\beta}}(w^k)+||\eta G^k||^2_2\nonumber\\
		\Longrightarrow&&\delta_{S\cap \overline{B}_{\beta}}(w^{k+1})+||w^{k+1}-w^k||^2_2+2\eta\langle w^{k+1}-w^k,G^k\rangle&\leq&&\delta_{S\cap \overline{B}_{\beta}}(w^k)\nonumber\\
		\Longrightarrow&&||w^{k+1}-w^k||^2_2+2\eta\langle w^{k+1}-w^k,G^k\rangle&\leq&&0\nonumber\\
		\Longrightarrow&&\frac{1}{2\eta}||w^{k+1}-w^k||^2_2+\langle w^{k+1}-w^k,G^k\rangle&\leq&&0,\label{proopt2}
	\end{alignat}
	where the third inequality holds since $w^k$ is feasible for $k\in\{1,...,K+1\}$. Continuing from \eqref{proopt2},
	\begin{alignat}{6}
		&f_{\alpha}(w^{k+1})+\frac{1}{2\eta}||w^{k+1}-w^k||^2_2+\langle w^{k+1}-w^k,G^k-\nabla f_{\alpha}(w^k)\rangle\nonumber\\
		\leq &f_{\alpha}(w^k)+\alpha^{-1}\sqrt{d}L_0||w^{k+1}-w^k||^2_2\label{xdiff}\\
		\Longrightarrow&2\eta^{-1}\langle w^{k+1}-w^k,G^k-\nabla f_{\alpha}(w^k)\rangle\nonumber\\
		\leq&2\eta^{-1}(f_{\alpha}(w^k)-f_{\alpha}(w^{k+1}))-\eta^{-2}||w^{k+1}-w^k||^2_2+2\alpha^{-1}\sqrt{d}L_0\eta^{-1}||w^{k+1}-w^k||^2_2,\label{join}
	\end{alignat}
	where the first inequality uses the smoothness of $f_{\alpha}(w)$ \citep[Lemma 1.2.3]{nesterov2004}:
	\begin{alignat}{6}
		f_{\alpha}(w^{k+1})\leq &f_{\alpha}(w^k)+\langle \nabla f_{\alpha}(w^k),w^{k+1}-w^k\rangle+\alpha^{-1}\sqrt{d}L_0||w^{k+1}-w^k||^2_2.\nonumber
	\end{alignat}
	Plugging \eqref{join} into \eqref{join0},
	\begin{alignat}{6}
		&&&\dist(0,\nabla f_{\alpha}(w^{k+1})+N(w^{k+1},S\cap \overline{B}_{\beta}))^2\nonumber\\
		&&\leq& (6\alpha^{-2}dL_0^2+4\alpha^{-1}\sqrt{d}L_0\eta^{-1}+\eta^{-2})||w^{k+1}-w^k||^2_2+3||\nabla f_{\alpha}(w^k)-G^k||^2_2\nonumber\\
		&&+&2\eta^{-1}(f_{\alpha}(w^k)-f_{\alpha}(w^{k+1}))-\eta^{-2}||w^{k+1}-w^k||^2_2+2\alpha^{-1}\sqrt{d}L_0\eta^{-1}||w^{k+1}-w^k||^2_2\nonumber\\
		&&=& (6\alpha^{-2}dL_0^2+6\alpha^{-1}\sqrt{d}L_0\eta^{-1})||w^{k+1}-w^k||^2_2+3||\nabla f_{\alpha}(w^k)-G^k||^2_2+2\eta^{-1}(f_{\alpha}(w^k)-f_{\alpha}(w^{k+1}))\nonumber\\
		&&=& 6\alpha^{-1}\sqrt{d}L_0(\alpha^{-1}\sqrt{d}L_0+\eta^{-1})||w^{k+1}-w^k||^2_2+3||\nabla f_{\alpha}(w^k)-G^k||^2_2+2\eta^{-1}(f_{\alpha}(w^k)-f_{\alpha}(w^{k+1})).\label{finalbound}
	\end{alignat}
	Rearranging \eqref{xdiff} and using Young's inequality again for the second inequality,
	\begin{alignat}{6}
		&&&(\frac{1}{2\eta}-\alpha^{-1}\sqrt{d}L_0)||w^{k+1}-w^k||^2_2\nonumber\\
		\leq&&&f_{\alpha}(w^k)-f_{\alpha}(w^{k+1})+\langle w^{k+1}-w^k,\nabla f_{\alpha}(w^k)-G^k\rangle\nonumber\\
		\leq&&&f_{\alpha}(w^k)-f_{\alpha}(w^{k+1})+\frac{\alpha^{-1}\sqrt{d}L_0}{2}||w^{k+1}-w^k||^2_2+\frac{\alpha}{2\sqrt{d}L_0}||\nabla f_{\alpha}(w^k)-G^k||^2_2\nonumber\\
		\Longrightarrow&&&(\frac{1}{2\eta}-\frac{3\alpha^{-1}\sqrt{d}L_0}{2})||w^{k+1}-w^k||^2_2\nonumber\\
		\leq&&&f_{\alpha}(w^k)-f_{\alpha}(w^{k+1})+\frac{\alpha}{2\sqrt{d}L_0}||\nabla f_{\alpha}(w^k)-G^k||^2_2\nonumber\\
		\Longrightarrow&&&\frac{1}{2}(\theta+\eta^{-1}-3\alpha^{-1}\sqrt{d}L_0)||w^{k+1}-w^k||^2_2\nonumber\\
		\leq&&&f_{\alpha}(w^k)-f_{\alpha}(w^{k+1})+\frac{\alpha}{2\sqrt{d}L_0}||\nabla f_{\alpha}(w^k)-G^k||^2_2+\frac{\theta}{2}||w^{k+1}-w^k||^2_2\label{wdiff}
	\end{alignat}
	for an arbitrary $\theta\in\RR$. Focusing on $||w^{k+1}-w^k||^2_2$,
	\begin{alignat}{6}
		&&&w^{k+1}\in\argmin_{x\in\RR^d}\text{ }\{\delta_{S\cap \overline{B}_{\beta}}(x)+||x-(w^k-\eta G^k)||_2\}\nonumber\\
		\Longrightarrow&&&||w^{k+1}-(w^k-\eta G^k)||_2\leq ||w^k-(w^k-\eta G^k)||_2=||\eta G^k||_2.\label{wdiffbound}
	\end{alignat}
	Using the reverse triangle inequality and \eqref{wdiffbound},
	\begin{alignat}{6}
		&&&||w^{k+1}-w^k||-||\eta G^k||_2\leq ||w^{k+1}-w^k+\eta G^k||_2\nonumber\\
		\Longrightarrow&&&||w^{k+1}-w^k||\leq 2||\eta G^k||_2,\nonumber
	\end{alignat}
	and applying this bound in \eqref{wdiff} assuming that $\theta\geq 0$,
	\begin{alignat}{6}
		&&&\frac{1}{2}(\theta+\eta^{-1}-3\alpha^{-1}\sqrt{d}L_0)||w^{k+1}-w^k||^2_2\nonumber\\
		\leq&&&f_{\alpha}(w^k)-f_{\alpha}(w^{k+1})+\frac{\alpha}{2\sqrt{d}L_0}||\nabla f_{\alpha}(w^k)-G^k||^2_2+2\theta\eta^2||G^k||^2_2.\label{wdiffbound2}
	\end{alignat}	
	
	Assuming that $\theta+\eta^{-1}-3\alpha^{-1}\sqrt{d}L_0>0$ and plugging \eqref{wdiffbound2} into \eqref{finalbound},
	\begin{alignat}{6}
		&&&\dist(0,\nabla f_{\alpha}(w^{k+1})+N(w^{k+1},S\cap \overline{B}_{\beta}))^2\nonumber\\
		&&\leq&  \frac{6\alpha^{-1}\sqrt{d}L_0(\alpha^{-1}\sqrt{d}L_0+\eta^{-1})}
		{\frac{1}{2}(\theta+\eta^{-1}-3\alpha^{-1}\sqrt{d}L_0)}\big(f_{\alpha}(w^k)-f_{\alpha}(w^{k+1})
		+\frac{\alpha}{2\sqrt{d}L_0}||\nabla f_{\alpha}(w^k)-G^k||^2_2\nonumber\\
		&&&+2\theta\eta^2||G^k||^2_2\big)+3||\nabla f_{\alpha}(w^k)-G^k||^2_2+2\eta^{-1}(f_{\alpha}(w^k)-f_{\alpha}(w^{k+1}))\nonumber\\
		\Longrightarrow&&&\EE[\dist(0,\nabla f_{\alpha}(w^{k+1})+N(w^{k+1},S\cap \overline{B}_{\beta}))^2]\nonumber\\
		&&\leq&  \frac{6\alpha^{-1}\sqrt{d}L_0(\alpha^{-1}\sqrt{d}L_0+\eta^{-1})}
		{\frac{1}{2}(\theta+\eta^{-1}-3\alpha^{-1}\sqrt{d}L_0)}\big(\EE[f_{\alpha}(w^k)-f_{\alpha}(w^{k+1})]
		+\frac{\alpha}{2\sqrt{d}L_0}\EE[||\nabla f_{\alpha}(w^k)-G^k||^2_2]\nonumber\\
		&&&+2\theta\eta^2\EE[||G^k||^2_2]\big)+3\EE[||\nabla f_{\alpha}(w^k)-G^k||^2_2]+2\eta^{-1}\EE[f_{\alpha}(w^k)-f_{\alpha}(w^{k+1})]\nonumber\\
		\Longrightarrow&&&\EE[\dist(0,\nabla f_{\alpha}(w^{R})+N(w^{R},S\cap \overline{B}_{\beta}))^2]\nonumber\\
		&&\leq&\frac{1}{K}\sum_{k=1}^{K}\bigg(  \frac{6\alpha^{-1}\sqrt{d}L_0(\alpha^{-1}\sqrt{d}L_0+\eta^{-1})}
		{\frac{1}{2}(\theta+\eta^{-1}-3\alpha^{-1}\sqrt{d}L_0)}
		\big(\EE[f_{\alpha}(w^k)-f_{\alpha}(w^{k+1})]
		+\frac{\alpha}{2\sqrt{d}L_0}\EE[||\nabla f_{\alpha}(w^k)-G^k||^2_2]\nonumber\\
		&&&+2\theta\eta^2\EE[||G^k||^2_2]\big)+3\EE[||\nabla f_{\alpha}(w^k)-G^k||^2_2]+2\eta^{-1}\EE[f_{\alpha}(w^k)-f_{\alpha}(w^{k+1})]\bigg)\nonumber\\
		&&\leq& \frac{6\alpha^{-1}\sqrt{d}L_0(\alpha^{-1}\sqrt{d}L_0+\eta^{-1})}
		{\frac{1}{2}(\theta+\eta^{-1}-3\alpha^{-1}\sqrt{d}L_0)}\bigg(\frac{\EE[f_{\alpha}(w^1)-
			f_{\alpha}(w^{K+1})]}{K}+\frac{\alpha}{2\sqrt{d}L_0}\frac{\upsilon Q}{M}+2\theta\eta^2\upsilon Q\bigg)\nonumber\\
		&&&+3\frac{\upsilon Q}{M}+\frac{2\eta^{-1}}{K}\EE[f_{\alpha}(w^1)-f_{\alpha}(w^{K+1})]\nonumber\\
		&&\leq& \frac{6\alpha^{-1}\sqrt{d}L_0(\alpha^{-1}\sqrt{d}L_0+\eta^{-1})}
		{\frac{1}{2}(\theta+\eta^{-1}-3\alpha^{-1}\sqrt{d}L_0)}\bigg(\frac{f_{\alpha}(w^1)-f_{\alpha}(w^*)}{K}
		+\big(\frac{\alpha}{2\sqrt{d}L_0M}+2\theta\eta^2\big)\upsilon Q\bigg)\nonumber\\
		&&&+3\frac{\upsilon Q}{M}+\frac{2\eta^{-1}}{K}(f_{\alpha}(w^1)-f_{\alpha}(w^*))\nonumber\\
		&&\leq& \frac{6\alpha^{-1}\sqrt{d}L_0(\alpha^{-1}\sqrt{d}L_0+\eta^{-1})}
		{\frac{1}{2}(\theta+\eta^{-1}-3\alpha^{-1}\sqrt{d}L_0)}\bigg(
		\frac{\Delta}{K}+\big(\frac{\alpha}{2\sqrt{d}L_0 M}+2\theta\eta^2\big)\upsilon Q\bigg)\nonumber\\
		&&&+3\frac{\upsilon Q}{M}+\frac{2\eta^{-1}\Delta}{K},\label{delta}
	\end{alignat}
	where the third last inequality uses Proposition \ref{dq} and the last inequality uses that $\Delta\geq f_{\alpha}(w^1)-f_{\alpha}(w^*)$. It holds that $f_{\alpha}(w^1)-f_{\alpha}(w^*)\leq 2\beta\sqrt{d}L_0$ from Proposition \ref{lipcont}.4 and 
	$f_{\alpha}(w^1)-f_{\alpha}(w^*_{\alpha})\leq f(w^1)-f_{\alpha}(w^*_{\alpha})+\alpha L_0\sqrt{\frac{d}{12}}\leq f(w^1)-f(w^*)+\alpha L_0\sqrt{\frac{d}{3}}$ using Propositions \ref{lipcont}.2 and \ref{lipcont}.3.\\
	
	Let $\theta=3\tau\alpha^{-1}\sqrt{d}L_0$ and $\eta^{-1}=3\rho\alpha^{-1}\sqrt{d}L_0$ for $\tau\geq 0$ and $\rho>0$ such that $\tau+\rho>1$. Continuing from \eqref{delta},
	\begin{alignat}{6}
		&&&\EE[\dist(0,\nabla f_{\alpha}(w^{R})+N(w^{R},S\cap \overline{B}_{\beta}))^2]\nonumber\\
		&&\leq& \frac{6\alpha^{-1}\sqrt{d}L_0(\alpha^{-1}\sqrt{d}L_0+3\rho\alpha^{-1}\sqrt{d}L_0)}
		{\frac{1}{2}(\tau+\rho-1)3\alpha^{-1}\sqrt{d}L_0}\bigg(\frac{\Delta}{K}+
		\big(\frac{\alpha}{2\sqrt{d}L_0 M}+\frac{\tau}{\rho^2}\frac{2\alpha}{3\sqrt{d}L_0}\big)\upsilon Q\bigg)\nonumber\\
		&&&+3\frac{\upsilon Q}{M}+\frac{6\rho\alpha^{-1}\sqrt{d}L_0\Delta}{K}\nonumber\\
		&&=& \frac{6\alpha^{-2}dL_0^2(1+3\rho)}
		{\frac{1}{2}(\tau+\rho-1)3\alpha^{-1}\sqrt{d}L_0}\bigg(\frac{\Delta}{K}+
		\big(\frac{1}{2M}+\frac{\tau}{\rho^2}\frac{2}{3}\big)\frac{\alpha}{\sqrt{d}L_0}\upsilon Q\bigg)+3\frac{\upsilon Q}{M}+\frac{6\rho\alpha^{-1}\sqrt{d}L_0\Delta}{K}\nonumber\\
		&&=& \frac{4\alpha^{-1}\sqrt{d}L_0(1+3\rho)}
		{(\tau+\rho-1)}\bigg(\frac{\Delta}{K}+\big(\frac{1}{2M}+\frac{\tau}{\rho^2}\frac{2}{3}\big)\frac{\alpha}{\sqrt{d}L_0}\upsilon Q\bigg)+3\frac{\upsilon Q}{M}+\frac{6\rho\alpha^{-1}\sqrt{d}L_0\Delta}{K}\nonumber\\
		&&=&\left(\frac{2(1+3\rho)}
		{(\tau+\rho-1)}+3\rho\right)\frac{2\alpha^{-1}\sqrt{d} L_0\Delta}{K}+\left(\frac{4(1+3\rho)}
		{(\tau+\rho-1)}\bigg(\frac{1}{2}+\frac{2}{3}\frac{M\tau}{\rho^2}\bigg)+3\right)\frac{\upsilon Q}{M}.\nonumber
	\end{alignat}
\end{proof}

\begin{proof}{\it (Corollary \ref{complexity})}
	In order for
	$|f_{\alpha}(w)-f(w)|\leq\epsilon_1\text{ for all }w\in \overline{B}_{\beta}$, we require $\alpha\leq\frac{\epsilon_1}{L_0\sqrt{\frac{d}{12}}}$ from Proposition \ref{lipcont}.2, and for $\widehat{\alpha}\leq \epsilon_3$, we require
	$\alpha\leq \frac{2\epsilon_3}{\sqrt{d}}$ from Corollary \ref{mainconeps}. To ensure that $\EE[\dist(0,\nabla f_{\alpha}(\overline{w})+N(\overline{w},S\cap \overline{B}_{\beta}))]\leq\epsilon$  or $\EE[\dist(0,\overline{\partial} f_{\widehat{\alpha}}(w^{R})+N(w^{R},S\cap \overline{B}_{\beta})]\leq\epsilon$ for $\epsilon=\epsilon_2$ or $\epsilon_4,$ using \eqref{theorem7} or \eqref{coroll8} and Jensen's inequality, it is sufficient for  
	\begin{alignat}{6}
		&&C_1\frac{2\alpha^{-1}\sqrt{d}L_0\Delta}{K}+C_2\frac{\upsilon Q}{M}\leq \epsilon^2.\nonumber
	\end{alignat}     
	Taking $y\in(0,1)$, we choose $K$ and $M$ such that  
	\begin{alignat}{6}
		&&C_1\frac{2\alpha^{-1}\sqrt{d}L_0\Delta}{K}\leq y\epsilon^2\nonumber
	\end{alignat}
	and
	\begin{alignat}{6}
		&&C_2\frac{\upsilon Q}{M}\leq (1-y)\epsilon^2,\nonumber
	\end{alignat}
	which results in requiring
	\begin{alignat}{6}
		&&C_1\frac{2\alpha^{-1}\sqrt{d}L_0\Delta}{y\epsilon^2}\leq K \label{K}
	\end{alignat}
	and
	\begin{alignat}{6}
		&&C_2\frac{\upsilon Q}{(1-y)\epsilon^2}\leq M.\nonumber
	\end{alignat}
	Assuming SPA is run for the full $K$ iterations, the number of gradient calls will equal $KM$. Considering the bound on KM, 
	\begin{alignat}{6}
		&&C_1\frac{2\alpha^{-1}\sqrt{d}L_0\Delta}{y\epsilon^2}C_2\frac{\upsilon Q}{(1-y)\epsilon^2}\leq KM,\label{KM}
	\end{alignat}
	minimizing the left-hand side of \eqref{KM} in terms of $y$ gives $y=0.5$, and minimizing the left-hand side of \eqref{K} in terms of $\alpha$ sets $\alpha=\frac{\epsilon_1}{L_0\sqrt{\frac{d}{12}}}$ or $\alpha= \frac{2\epsilon_3}{\sqrt{d}}$. Using these values for $\alpha$ and $y$ gives the bounds for $K$ of
	\begin{alignat}{6}
		&&C_1\sqrt{\frac{4}{3}}\frac{d L_0^2\Delta}{\epsilon_1\epsilon^2_2}\leq K\nonumber
	\end{alignat}
	for an $(\epsilon_1,\epsilon_2)$-solution, and
	\begin{alignat}{6}
		&&C_1\frac{2dL_0\Delta}{\epsilon_3\epsilon^2_4}\leq K\nonumber
	\end{alignat}
	for an $(\epsilon_3,\epsilon_4)$-solution. The bound for $M$ equals
	\begin{alignat}{6}
		&&C_2\frac{2\upsilon Q}{\epsilon^2}\leq M\nonumber
	\end{alignat}
	for $\epsilon=\epsilon_2$ or $\epsilon_4$.\\
	
	Taking the choices of $K$ and $M$ given in this corollary and the upper bound of $KM$ for the total number of gradient calls gives the gradient call complexities of $O(\epsilon_1^{-1}\epsilon_2^{-4})$ and $O(\epsilon_3^{-1}\epsilon_4^{-4})$ to achieve an expected $(\epsilon_1,\epsilon_2)$ and $(\epsilon_3,\epsilon_4)$-stationary point, respectively. Given that one projection is done per iteration, the projection operator complexities are  $O(\epsilon_1^{-1}\epsilon_2^{-2})$ and $O(\epsilon_3^{-1}\epsilon_4^{-2})$.
\end{proof}

\begin{proof}{\it (Corollary \ref{comp2})}
	Following \citep[Equation 2.28]{ghadimi2013},
	\begin{alignat}{6}
		&&&\dist(0,G(w^*)+N(w^*,S\cap \overline{B}_{\beta}))^2\nonumber\\
		&&=\min_{w\in W}&\dist(0,G(w)+N(w,S\cap \overline{B}_{\beta}))^2\nonumber\\
		&&=\min_{w\in W}&\dist(0,\nabla f_{\alpha}(w)+N(w,S\cap \overline{B}_{\beta})+G(w)-\nabla f_{\alpha}(w))^2\nonumber\\
		&&\leq\min_{w\in W}&\{2\dist(0,\nabla f_{\alpha}(w)+N(w,S\cap \overline{B}_{\beta}))^2+2||G(w)-\nabla f_{\alpha}(w)||^2_2\}\nonumber\\
		&&\leq\min_{w\in W}&2\dist(0,\nabla f_{\alpha}(w)+N(w,S\cap \overline{B}_{\beta}))^2+\max_{w\in W}2||G(w)-\nabla f_{\alpha}(w)||^2_2,\label{hp1}
	\end{alignat}
	where the first inequality holds since
	
	\begin{alignat}{6}
		&&&\dist(0,\nabla f_{\alpha}(w)+N(w,S\cap \overline{B}_{\beta})+G(w)-\nabla f_{\alpha}(w))^2\nonumber\\
		&&=\min_{\nu\in N(w,S\cap\overline{B}_{\beta})}&||\nabla f_{\alpha}(w)+\nu+G(w)-\nabla f_{\alpha}(w)||^2_2\nonumber\\
		&&=\min_{\nu\in N(w,S\cap\overline{B}_{\beta})}&||\nabla f_{\alpha}(w)+\nu||^2_2+2\langle \nabla f_{\alpha}(w)+\nu,G(w)-\nabla f_{\alpha}(w)\rangle+||G(w)-\nabla f_{\alpha}(w)||^2_2\nonumber\\
		&&\leq\min_{\nu\in N(w,S\cap\overline{B}_{\beta})}&2||\nabla f_{\alpha}(w)+\nu||^2_2+2||G(w)-\nabla f_{\alpha}(w)||^2_2\nonumber
	\end{alignat}
	using Young's inequality. Using Young's inequality again for the first inequality and \eqref{hp1} for the second inequality, 
	
	\begin{alignat}{6}
		&&&\dist(0,\nabla f_{\alpha}(w^*)+N(w^*,S\cap \overline{B}_{\beta}))^2\nonumber\\
		&&\leq2&\dist(0,G(w^*)+N(w^*,S\cap \overline{B}_{\beta}))^2+2||\nabla f_{\alpha}(w^*)-G(w^*)||^2_2\nonumber\\
		&&\leq4&\min_{w\in W}\dist(0,\nabla f_{\alpha}(w)+N(w,S\cap \overline{B}_{\beta}))^2+4\max_{w\in W}||G(w)-\nabla f_{\alpha}(w)||^2_2+2||\nabla f_{\alpha}(w^*)-G(w^*)||^2_2\nonumber\\
		&&\leq4&\min_{w\in W}\dist(0,\nabla f_{\alpha}(w)+N(w,S\cap \overline{B}_{\beta}))^2+6\max_{w\in W}||G(w)-\nabla f_{\alpha}(w)||^2_2.\label{hp3}
	\end{alignat}
	
	Let the right-hand side of \eqref{theorem7} be denoted as 
	\begin{alignat}{6}
		D:= C_1\frac{2\alpha^{-1}\sqrt{d} L_0\Delta}{K}+C_2\frac{\upsilon Q}{M}.\nonumber
	\end{alignat}
	
	Considering the first term of \eqref{hp3}, 
	\begin{alignat}{6}
		&&&\PP\bigg(4\min_{w\in W}\dist(0,\nabla f_{\alpha}(w)+N(w,S\cap \overline{B}_{\beta}))^2\geq 4e D\bigg)\nonumber\\
		&&=&\Pi_{i=1}^{r}\PP(4\dist(0,\nabla f_{\alpha}(w^i)+N(w^i,S\cap \overline{B}_{\beta}))^2\geq4e D)\nonumber\\
		&&\leq& e^{-r}\nonumber
	\end{alignat}
	using Markov's inequality. For the second term of \eqref{hp3}, using Proposition \ref{dq} and Boole's inequality,
	\begin{alignat}{6}
		&&&\PP\bigg(6\max_{w\in W}||G(w)-\nabla f_{\alpha}(w)||^2_2\geq 6\psi \frac{\upsilon Q}{T}\bigg)\nonumber\\
		=&&&\PP\bigg(\bigcup_{i=1}^r \bigg\{6||G(w^i)-\nabla f_{\alpha}(w^i)||^2_2\geq 6\psi \frac{\upsilon Q}{T}\bigg\}\bigg)\nonumber\\
		\leq&&&\frac{r}{\psi}.\nonumber
	\end{alignat}
	Combining these two probability inequalities together, for the left-hand side of \eqref{hp3}, 
	\begin{alignat}{6}
		&&&\PP\bigg(\dist(0,\nabla f_{\alpha}(w^*)+N(w^*,S\cap \overline{B}_{\beta}))^2\geq 4eD+6\psi\frac{\upsilon Q}{T}\bigg)\nonumber\\
		&&\leq&\PP\bigg(
		4\min_{w\in W}\dist(0,\nabla f_{\alpha}(w)+N(w,S\cap \overline{B}_{\beta}))^2+6\max_{w\in W}||G(w)-\nabla f_{\alpha}(w)||^2_2\geq 4e D+6\psi\frac{\upsilon Q}{T}\bigg)\nonumber\\
		&&\leq&\PP\bigg(\{4\min_{w\in W}\dist(0,\nabla f_{\alpha}(w)+N(w,S\cap \overline{B}_{\beta}))^2\geq 4e D\}\nonumber\\
		&&&\cup\{6\max_{w\in W}||G(w)-\nabla f_{\alpha}(w)||^2_2\geq6\psi\frac{\upsilon Q}{T}\}\bigg)\nonumber\\
		&&\leq& e^{-r}+\frac{r}{\psi}.\label{prob_bound}
	\end{alignat}
	
	An upper bound on the total number of gradient calls required for computing $W$ and $G(w)$ for $w\in W$ is equal to $r(KM+T)$. Using \eqref{prob_bound}, the minimization of $r(KM+T)$ while ensuring that $\PP(\dist(0,\nabla f_{\alpha}(w^*)+N(w^*,S\cap \overline{B}_{\beta}))> \epsilon_{2/4})\leq \gamma$, where $\alpha=\frac{\epsilon_1}{L_0\sqrt{\frac{d}{12}}}$ and $\epsilon_{2/4}=\epsilon_2$, or $\alpha=\frac{2\epsilon_3}{\sqrt{d}}$ and $\epsilon_{2/4}=\epsilon_4$ (assuming that $\kappa>\beta+\epsilon_3$) can be written as
	\begin{alignat}{6}
		\min\limits_{\substack{r,K,M,\\T,\psi}}&\text{ }&&r(KM+T)\nonumber\\
		\st&&&4e \bigg(C_1\frac{2\alpha^{-1}\sqrt{d} L_0\Delta}{K}+C_2\frac{\upsilon Q}{M}\bigg)+6\psi\frac{\upsilon Q}{T}\leq \epsilon^2_{2/4}\label{eps2eq}\\
		&&&e^{-r}+\frac{r}{\psi}\leq \gamma\nonumber\\
		&&& r,K,M,T\in\ZZ_{>0},\quad\psi>0.\nonumber 
	\end{alignat}
	
	Writing \eqref{eps2eq} as
	\begin{alignat}{6}
		&&&C_1\frac{2\alpha^{-1}\sqrt{d} L_0\Delta}{K}+C_2\frac{\upsilon Q}{M}\leq \frac{\epsilon^2_{2/4}-6\psi\frac{\upsilon Q}{T}}{4e},\nonumber
	\end{alignat}
	
	we set $K$ and $M$ according to Corollary \ref{complexity} to find an expected $(\epsilon_1,\epsilon_2')$ or $(\epsilon_3,\epsilon_4')$-stationary point for $\epsilon_2'=\sqrt{\frac{\epsilon^2_2-6\psi\frac{\upsilon Q}{T}}{4e}}$ or 
	$\epsilon_4'=\sqrt{\frac{\epsilon^2_4-6\psi\frac{\upsilon Q}{T}}{4e}}$
	assuming that $\epsilon^2_{2/4}-6\psi\frac{\upsilon Q}{T}>0$:
	\begin{alignat}{6}
		&&K^*=\left\lceil C_1\frac{2}{\chi}\frac{d L_0^2\Delta}{\epsilon_{1/3}(\epsilon_{2/4}')^2}\right\rceil\quad\text{and}\quad M^*=\left\lceil C_2\frac{2\upsilon Q}{(\epsilon_{2/4}')^2}\right\rceil,\nonumber
	\end{alignat}
	where $\epsilon_{1/3}=\epsilon_1$ or $\epsilon_3$, $\epsilon_{2/4}'=\epsilon_2'$ or $\epsilon_4'$, and $\chi=\sqrt{3}$ or $L_0$ for an expected  $(\epsilon_1,\epsilon_2')$ or $(\epsilon_3,\epsilon_4')$-stationary point, respectively. The optimization problem then becomes
	\begin{alignat}{6}
		\min\limits_{r,T,\psi}&\text{ }&&r(K^*M^*+T)\nonumber\\
		\st&&&e^{-r}+\frac{r}{\psi}\leq \gamma\label{rphi} \\
		&&&\epsilon^2_{2/4}-6\psi\frac{\upsilon Q}{T}>0\label{Tineq}\\
		&&&r,T\in\ZZ_{>0},\quad \psi>0.\nonumber 
	\end{alignat}
	
	For any $c\in(0,1)$, let 
	\begin{alignat}{6}
		e^{-r}\leq c\gamma\quad\text{and}\quad\frac{r}{\psi}\leq (1-c)\gamma.\nonumber
	\end{alignat}
	Setting $r$ and $\psi$ to $r^*=\lceil-\ln(c\gamma)\rceil$ and $\psi^*=\frac{\lceil-\ln(c\gamma)\rceil}{(1-c)\gamma}$ is then valid for \eqref{rphi}. For any $\phi>1$, setting $T$ to $T^*=\lceil6\phi\psi^*\frac{\upsilon Q}{\epsilon^2_{2/4}}\rceil$ is valid for \eqref{Tineq}. The total number of gradient calls then equals
	\begin{alignat}{6}
		&&&r^*\bigg(K^*M^*+T^*\bigg)\nonumber\\
		&&=&r^*\bigg(\left\lceil C_1\frac{2}{\chi}\frac{d L_0^2\Delta}{\epsilon_{1/3}(\epsilon_{2/4}')^2}\right\rceil\left\lceil C_2\frac{2\upsilon Q}{(\epsilon_{2/4}')^2}\right\rceil+T^*\bigg)\nonumber\\
		&&=&r^*\bigg(\left\lceil C_1\frac{2}{\chi}\frac{d L_0^2\Delta}{\epsilon_{1/3}\frac{\epsilon^2_{2/4}-6\psi^*\frac{\upsilon Q}{T^*}}{4e}}\right\rceil\left\lceil C_2\frac{2\upsilon Q}{\frac{\epsilon^2_{2/4}-6\psi^*\frac{\upsilon Q}{T^*}}{4e}}\right\rceil+T^*\bigg)\nonumber\\
		&&\leq&r^*\bigg(\left\lceil C_1\frac{2}{\chi}\frac{4ed L_0^2\Delta}{\epsilon_{1/3}\epsilon^2_{2/4}(1-\phi^{-1})}\right\rceil\left\lceil C_2\frac{8e\upsilon Q}{\epsilon^2_{2/4}(1-\phi^{-1})}\right\rceil+T^*\bigg)\nonumber\\
		&&=&\lceil-\ln(c\gamma)\rceil\bigg(\left\lceil C_1\frac{2}{\chi}\frac{4ed L_0^2\Delta}{\epsilon_{1/3}\epsilon^2_{2/4}(1-\phi^{-1})}\right\rceil\left\lceil C_2\frac{8e\upsilon Q}{\epsilon^2_{2/4}(1-\phi^{-1})}\right\rceil+\bigg\lceil6\phi\frac{\lceil-\ln(c\gamma)\rceil}{(1-c)\gamma}\frac{\upsilon Q}{\epsilon^2_{2/4}}\bigg\rceil\bigg), \nonumber
	\end{alignat}
	where the inequality holds since 
	\begin{alignat}{6}
		&&&\epsilon^2_{2/4}-6\psi^*\frac{\upsilon Q}{T^*}=\epsilon^2_{2/4}-6\psi^*\frac{\upsilon Q}{\lceil6\phi\psi^*\frac{\upsilon Q}{\epsilon^2_{2/4}}\rceil}\geq\epsilon^2_{2/4}-6\psi^*\frac{\upsilon Q}{6\phi\psi^*\frac{\upsilon Q}{\epsilon^2_{2/4}}}=\epsilon^2_{2/4}(1-\phi^{-1}),\nonumber
	\end{alignat}
	and the gradient call complexity equals
	$\tilde{O}\left(\epsilon^{-1}_1\epsilon^{-4}_2+\gamma^{-1}\epsilon^{-2}_2\right)$ or 
	$\tilde{O}\left(\epsilon^{-1}_3\epsilon^{-4}_4+\gamma^{-1}\epsilon^{-2}_4\right)$ for an $(\epsilon_1,\epsilon_2)$ or $(\epsilon_3,\epsilon_4)$-stationary point with probability $1-\gamma$. The number of projection computations is upper bounded by $r^*K^*$ which has a complexity of $\tilde{O}\left(\epsilon^{-1}_1\epsilon^{-2}_2\right)$ or   
	$\tilde{O}\left(\epsilon^{-1}_3\epsilon^{-2}_4\right)$. 
\end{proof}

\begin{proof}{\it (Proposition \ref{normcomp})}	
	For simplicity let $I:= I(w)$. The distance function 
	\begin{alignat}{6}
		\dist(0,G(w)+N(w,C\cap \overline{B}_{\beta}))=\min_{v\in N(w,C\cap \overline{B}_{\beta})}\text{ }&&&||G(w)+v||_2.\nonumber
	\end{alignat}	
	When $Y=I$, we only need to focus on two subsets of $[n]$: $I$ and $[n]\setminus I$. Each subset of elements $v^i$ for $i\in I$ is set to its optimal value while following Proposition \ref{limitnorm}: when $w^i_j=\beta$, it is required that $v^i_j\geq 0$, hence it is optimal to set $v^i_j=0$ when $G^i_j(w)\geq 0$ and to set $v^i_j=-G^i_j(w)$ when $G^i_j(w)<0$. Similarly, when $w^i_j=-\beta$, it is required that $v^i_j\leq 0$, and it is optimal to set $v^i_j=0$ when $G^i_j(w)\leq 0$ and to set $v^i_j=-G^i_j(w)$ when $G^i_j(w)>0$. When $i\in I$ and $|w^i_j|<\beta$, it is required that $v^i_j=0$. When $i\notin I$, $v^i_j$ is a free variable which is optimally set to $v^i_j=-G^i_j(w)$ for all $j\in[d_i]$.\\
	
	When there exists an $X\in Y$ such that $X\setminus I(w)\neq\{\emptyset\}$, and assuming that $\{p_i\}\subset\QQ_{>0}$, $C$ can be written equivalently with parameters $\{p'_i\}=\{c p_i\}$ and $m'=c m$ for a $c\in\ZZ_{>0}$ sufficiently large such that $\{p'_i\}\subset \ZZ_{>0}$, so without loss of generality we can assume that $\{p_i\}\subset \ZZ_{>0}$.\\
	
	The distance function has the same optimal solutions as 
	\begin{alignat}{6}
		\min_{v\in N(w,C\cap \overline{B}_{\beta})}\text{ }&&&||G(w)+v||^2_2,\nonumber
	\end{alignat}
	which can be written as 
	\begin{alignat}{6}
		\min\text{ }&&&||G(w)+v||^2_2\nonumber\\
		\text{s.t. }&&&\sum_{i\in I}p_i+\sum_{i\notin I}z_ip_i\leq m\nonumber\\	
		&&&\sum_{i\in I}p_i+\sum_{i\notin I}z_ip_i+(1-z_j)p_j\geq (1-z_j)(\lfloor m\rfloor+1)\quad \forall j\notin I\nonumber\\
		&&&v^i_j=\begin{cases}
			-G^i_j(w) & \text{ if }U^i_j\\ 
			0 & \text{otherwise}
		\end{cases}\quad \forall i\in I, \forall j\in[d_i]\nonumber\\	
		&&&v^i_j=-G^i_j(w)(1-z_i)\quad \forall i\notin I, \forall j\in[d_i]\nonumber\\
		&&&z_i\in\{0,1\}\text{ }\forall i\notin I.\nonumber
	\end{alignat}
	
	Using Proposition \ref{limitnorm}, the first two constraints determine an $X$ equal to $I$ and the indices for which $z_i=1$, where the second constraint enforces that $\sum_{i\in I}p_i+\sum_{i\notin I}z_ip_i+p_j>m$, given the integrality of $\{p_i\}$, for $j\notin X$, i.e. $z_j=0$. The third and fourth constraints set each subset of elements $v^i$ for $i\in I$ and $i\notin I$ to their optimal value, respectively. For the fourth constraint, when $z_i=1$, i.e. $i\in X$, $v^i$ must be set to $v^i=0$ given that $w^i=0$, and when $z_i=0$, $v^i$ is free to be chosen as $v^i=-G^i(w)$ to minimize the objective.\\ 
	
	The final binary integer program to compute $\dist(0,G(w)+N(w,C\cap \overline{B}_{\beta})$ is as follows, removing the $v$ decision variables.
	
	\begin{alignat}{6}
		\min\text{ }&&&\sum_{i\notin I}||G^i(w)||^2_2z_i\label{fin_binp}\\
		\text{s.t. }&&&\sum_{i\in I}p_i+\sum_{i\notin I}z_ip_i\leq m\nonumber\\	
		&&&\sum_{i\in I}p_i+\sum_{i\notin I}z_ip_i+(1-z_j)p_j\geq (1-z_j)(\lfloor m\rfloor+1)\quad \forall j\notin I\nonumber\\
		&&&z_i\in\{0,1\}\text{ }\forall i\notin I.\nonumber
	\end{alignat}
	Given an optimal solution $z^*$ to \eqref{fin_binp}, let $y^*$ be defined as
	\begin{alignat}{6}
		&&&(y^i_j)^*=\begin{cases}
			0 & \text{ if }U^i_j\\ 
			1 & \text{otherwise}
		\end{cases}\quad \forall i\in I, \forall j\in[d_i]\nonumber\\	
		&&&(y^i_j)^*=z^*_i\quad \forall i\notin I, \forall j\in[d_i].\nonumber
	\end{alignat}
	It follows that $\dist(0,G(w)+N(w,C\cap \overline{B}_{\beta})=\sqrt{\sum_{i\in[n]}\sum_{j\in d_i}(G^i_j(w))^2(y^i_j)^*}$.
\end{proof}

\begin{proof}{\it (Theorem \ref{asymptotic})}
	If $\beta$ is finite, $\{w^i\}\subset \overline{B}_{\beta}$ is a bounded sequence and there exists an accumulation point $\overline{w}$ of $\{w^i\}$. Otherwise, assume there exists an accumulation point $\overline{w}$ of $\{w^i\}$. For simplicity, let $\{w^i\}$ be redefined as a subsequence of $\{w^i\}$ such that $\lim\limits_{i\rightarrow\infty}w^i=\overline{w}$. Since $||\zeta||_2\leq L_0$ for all $\zeta\in \overline{\partial} f(w)$ for all $w\in \overline{B}_{\beta}$ \citep[Proposition 2.1.2]{clarke1990} and $0\in N(w,S\cap \overline{B}_{\beta})$,
	\begin{alignat}{6}
		\dist(0,\overline{\partial} f(w)+N(w,S\cap \overline{B}_{\beta}))=\dist(0,\overline{\partial} f(w)+N(w,S\cap \overline{B}_{\beta})\cap \overline{B}(0,2L_0))\nonumber
	\end{alignat} 
	for $w\in \overline{B}_{\beta}$: If $\zeta\in \overline{\partial} f(w)$, $\nu\in N(w,S\cap \overline{B}_{\beta})$, and $||\nu||_2>2L_0$, then by the reverse triangle inequality, $||\zeta+0||_2\leq L_0<||\nu||_2-||\zeta||_2\leq ||\nu+\zeta||_2$.\\
	
	Given that $N(w,S\cap \overline{B}_{\beta})$ is outer semicontinuous (proof of Proposition \ref{meas}), for all $\omega_1>0$, there exists an $\omega_2>0$ such that 
	\begin{alignat}{6}
		N(\hat{w},S\cap \overline{B}_{\beta})\cap \overline{B}(0,2L_0)\subseteq N(w,S\cap \overline{B}_{\beta})+\overline{B}(0,\omega_1)\label{cusp}
	\end{alignat}
	for all $\hat{w}\in \overline{B}(w,\omega_2)$ \citep[Proposition 5.12]{rockafellar2009}.\\ 
	
	For any $i\in \NN$ there exists an $I\in \NN$ such that for all $j>I$, $w^j$ is an expected $\left(\frac{\epsilon^i_3}{2},\frac{\epsilon^i_4}{2i}\right)$-stationary point and $||w^j-\overline{w}||_2\leq \min(\frac{\epsilon^i_3}{2},\omega^i_2)$, where $\omega^i_2>0$ is chosen such that \eqref{cusp} holds for $\omega_1=\frac{\epsilon^i_4}{2i}$, $\omega_2=\omega^i_2$, and $w=\overline{w}$. For such $w^j$, $\overline{B}(w^j,\frac{\epsilon^i_3}{2})\subset \overline{B}(\overline{w},\epsilon^i_3)$, hence 
	\begin{alignat}{6}  
		\{\overline{\partial} f(w): w\in \overline{B}(w^j,\epsilon^i_3/2)\}\subseteq\{\overline{\partial} f(w): w\in \overline{B}(\overline{w},\epsilon^i_3)\},\nonumber
	\end{alignat}
	and $\overline{\partial}_{\epsilon^i_3/2} f(w^j)\subseteq \overline{\partial}_{\epsilon^i_3} f(\overline{w})$. In addition,  
	\begin{alignat}{6}
		\frac{\epsilon^i_4}{2i}&\geq&&\EE[\dist(0,\overline{\partial}_{\epsilon^i_3/2} f(w^j)+N(w^j,S\cap \overline{B}_{\beta}))]\nonumber\\
		&=&&\EE[\dist(0,\overline{\partial}_{\epsilon^i_3/2} f(w^j)+N(w^j,S\cap \overline{B}_{\beta})\cap \overline{B}(0,2L_0))]\nonumber\\
		&\geq&&\EE[\dist(0,\overline{\partial}_{\epsilon^i_3}f(\overline{w})+N(w^j,S\cap \overline{B}_{\beta})\cap \overline{B}(0,2L_0))]\nonumber\\
		&\geq&&\EE[\dist(0,\overline{\partial}_{\epsilon^i_3}f(\overline{w})
		+N(\overline{w},S\cap \overline{B}_{\beta})+\overline{B}(0,\epsilon^i_4/(2i)))].\label{stabound}  
	\end{alignat}
	Using the reverse triangle inequality for the first inequality,
	\begin{alignat}{6}
		&&&\dist(0,\overline{\partial}_{\epsilon^i_3}f(\overline{w})
		+N(\overline{w},S\cap \overline{B}_{\beta})+\overline{B}(0,\epsilon^i_4/(2i)))\nonumber\\
		=&&&\min_{\substack{z\in\overline{\partial}_{\epsilon^i_3}f(\overline{w})+N(\overline{w},S\cap \overline{B}_{\beta})\\y\in \overline{B}(0,\epsilon^i_4/(2i))}}||z+y||_2\nonumber\\  
		\geq&&&\min_{\substack{z\in\overline{\partial}_{\epsilon^i_3}f(\overline{w})+N(\overline{w},S\cap \overline{B}_{\beta})\\y\in \overline{B}(0,\epsilon^i_4/(2i))}}||z||_2-||y||_2\nonumber\\ 
		=&&&\dist(0,\overline{\partial}_{\epsilon^i_3}f(\overline{w})
		+N(\overline{w},S\cap \overline{B}_{\beta}))-\frac{\epsilon^i_4}{2i}\nonumber\\ 
		\Longrightarrow&&&\EE[\dist(0,\overline{\partial}_{\epsilon^i_3}f(\overline{w})
		+N(\overline{w},S\cap \overline{B}_{\beta})+\overline{B}(0,\epsilon^i_4/(2i)))]\nonumber\\
		\geq&&&\EE[\dist(0,\overline{\partial}_{\epsilon^i_3}f(\overline{w})
		+N(\overline{w},S\cap \overline{B}_{\beta}))]-\frac{\epsilon^i_4}{2i}.\label{stabound2}
	\end{alignat}
	Applying \eqref{stabound2} in \eqref{stabound},
	\begin{alignat}{6}
		&&&\EE[\dist(0,\overline{\partial}_{\epsilon^i_3}f(\overline{w})
		+N(\overline{w},S\cap \overline{B}_{\beta}))]\leq \frac{\epsilon^i_4}{i}.\nonumber
	\end{alignat}
	From Markov's inequality,
	\begin{alignat}{6}
		\PP[\dist(0,\overline{\partial}_{\epsilon^i_3}f(\overline{w})+N(\overline{w},S\cap \overline{B}_{\beta})\geq\frac{1}{i}]\leq \epsilon^i_4.\nonumber
	\end{alignat}
	The sets 
	\begin{alignat}{6}
		V_i:=\{\overline{w}\in\RR^d:\dist(0,\overline{\partial}_{\epsilon^i_3}f(\overline{w})+N(\overline{w},S\cap \overline{B}_{\beta})\geq \frac{1}{i}\}\nonumber
	\end{alignat}
	are monotonically increasing: $V_i\subseteq V_{i+1}$, 
	as $\dist(0,\overline{\partial}_{\epsilon^i_3}f(\overline{w})+N(\overline{w},S\cap \overline{B}_{\beta})\leq \dist(0,\overline{\partial}_{\epsilon^{i+1}_3}f(\overline{w})+N(\overline{w},S\cap \overline{B}_{\beta})$ and $\frac{1}{i}>\frac{1}{i+1}$. The limit $\lim\limits_{i\rightarrow \infty}V_i=\bigcup\limits_{i\geq 1}V_i$ exists \citep[Excecise 2.F.]{bartle1995}, and is Borel measurable as a countable union of measurable sets, given that the functions $\dist(0,\overline{\partial}_{\epsilon^i_3}f(\overline{w})+N(\overline{w},S\cap \overline{B}_{\beta})$ are Borel measurable from Proposition \ref{meas}.\\ 
	
	We now want to prove that $\lim\limits_{i\rightarrow \infty}V_i=\{\overline{w}\in\RR^d:\dist(0,\overline{\partial}f(\overline{w})+N(\overline{w},S\cap \overline{B}_{\beta})>0\}$. For any $w\in \bigcup\limits_{i\geq 1}V_i$ there exists an $i\geq 1$ such that $\dist(0,\overline{\partial}f(w)+N(w,S\cap \overline{B}_{\beta})\geq \dist(0,\overline{\partial}_{\epsilon^i_3}f(w)+N(w,S\cap \overline{B}_{\beta})\geq \frac{1}{i}>0$, hence $w\in \{\overline{w}\in\RR^d:\dist(0,\overline{\partial}f(\overline{w})+N(\overline{w},S\cap \overline{B}_{\beta})>0\}$.\\ 
	
	For a $w\in \{\overline{w}\in\RR^d:\dist(0,\overline{\partial}f(\overline{w})+N(\overline{w},S\cap \overline{B}_{\beta})>0\}$, let $\omega=\dist(0,\overline{\partial}f(w)+N(w,S\cap \overline{B}_{\beta})$. As was shown with $N(w,S\cap \overline{B}_{\beta})$, given that the Clarke subdifferential is an upper semicontinuous set valued mapping \citep[Proposition 2.1.5 (d)]{clarke1990}, for all $\omega_1>0$, there exists an $\omega_2>0$ such that $\overline{\partial} f(\hat{w})\subset \overline{\partial} f(w)+\overline{B}(0,\omega_1)$ for all $\hat{w}\in \overline B(w,\omega_2)$, from which it follows that $\overline{\partial}_{\omega_2} f(w)\subseteq \text{co}\{\overline{\partial} f(w)+\overline{B}(0,\omega_1)\}=\overline{\partial} f(w)+\overline{B}(0,\omega_1)$, given that $\overline{\partial} f(w)$ is convex, and 
	\begin{alignat}{6}
		\dist(0,\overline{\partial}_{\omega_2} f(w)+N(w,S\cap \overline{B}_{\beta}))\geq \dist(0,\overline{\partial} f(w)+\overline{B}(0,\omega_1)+N(w,S\cap \overline{B}_{\beta})).\label{omega2dist}
	\end{alignat}
	Just as in proving \eqref{stabound2}, $\dist(0,\overline{\partial} f(w)+\overline{B}(0,\omega_1)+N(w,S\cap \overline{B}_{\beta}))$ can be bounded below: 
	\begin{alignat}{6}
		\dist(0,\overline{\partial} f(w)+\overline{B}(0,\omega_1)+N(w,S\cap \overline{B}_{\beta}))\geq\omega -\omega_1.\label{lowerbound}
	\end{alignat}
	Choosing $\omega_1=\frac{\omega}{2}$, there exists an $\omega_2>0$ such that  
	$\dist(0,\overline{\partial}_{\omega_2} f(w)+N(w,S\cap \overline{B}_{\beta}))\geq \frac{\omega}{2}$ from \eqref{omega2dist} and \eqref{lowerbound}. A $J\in \NN$ exists such that for all $i\geq J$, $\epsilon^i_3\leq\omega_2$. Setting $I\geq\max\{J,\frac{2}{\omega}\}$, 
	\begin{alignat}{6}
		\dist(0,\overline{\partial}_{\epsilon^i_3} f(w)+N(w,S\cap \overline{B}_{\beta}))\geq \frac{1}{i},\nonumber
	\end{alignat}
	
	i.e. $w\in V_i$, for all $i\geq I$, proving that $\lim\limits_{i\rightarrow \infty}V_i=\{\overline{w}\in\RR^d:\dist(0,\overline{\partial}f(\overline{w})+N(\overline{w},S\cap \overline{B}_{\beta})>0\}$.\\ 
	
	It holds that $\PP[\dist(0,\overline{\partial} f(\overline{w})+N(\overline{w},S\cap \overline{B}_{\beta}))=0]=1$ as 
	\begin{alignat}{6}
		&&\PP[\dist(0,\overline{\partial} f(\overline{w})+N(\overline{w},S\cap \overline{B}_{\beta}))>0]&=\lim\limits_{i\rightarrow\infty}\PP[\dist(0,\overline{\partial}_{\epsilon^i_3}f(\overline{w})+N(\overline{w},S\cap \overline{B}_{\beta})\geq \frac{1}{i}]\nonumber\\
		&&&\leq\lim_{i\rightarrow\infty}\epsilon^i_4=0,\nonumber
	\end{alignat}
	where the equality holds given that $V_i\subseteq V_{i+1}$ \citep[Theorem A.1.1]{shreve2004}.
\end{proof}

\section*{Appendix E. Section \ref{backprop} Proofs}

\begin{proof}{\it (Proposition \ref{boltebackprop})}
	For all $j\in I_k$, $\overline{\partial}\sigma_{j}(\cdot)$ is outer semicontinuous (\citealp[Proposition 2.1.5 (d)]{clarke1990}; \citealp[Theorem 5.19]{rockafellar2009}), hence measurable \citep[Exercise 14.9]{rockafellar2009}. There then exists a measurable selection $\widetilde{\nabla}\sigma_{j}(\cdot)\in \overline{\partial}\sigma_{j}(\cdot)$ for all $j\in I_k$ \citep[Theorem 14.6]{rockafellar2009}.\\ 
	
	As the product and sum of real-valued Borel measurable functions, see \citep[Algorithm 3]{bolte2021}, $\widetilde{\nabla}F(w,\xi_k)$ is Borel measurable in $w\in\RR^d$ for each $\xi_k\in\{\xi_i\}_{i=1}^{\infty}$.
	By the assumption that for each $i\in\{1,2,...\}$, all $\{\sigma_{j}\}_{j\in I_i}$ are definable in the same o-minimal structure, the Clarke subdifferentials $\{\overline{\partial}\sigma_{j}\}_{j\in I_i}$ are definable conservative fields in said o-minimal structure as well \citep[Remark 8]{bolte2021}, hence for each $\xi_k\in\{\xi_i\}_{i=1}^{\infty}$, $\widetilde{\nabla}F(w,\xi_k)$ will equal the gradient of $F(w,\xi_k)$ for almost every $w\in\RR^d$ following \citep[Corollary 5]{bolte2021}. Further, $\widetilde{\nabla}F(w,\xi)$ will equal the gradient of $F(w,\xi)$ for all $(w,\xi)\in\RR^{d+p}$ except 
	on a countable number of null sets: The set $\{(w,\xi): w\in\RR^d, \xi\notin\{\xi_i\}_{i=1}^{\infty}\}$ 
	which has measure zero by assumption, and potentially a null set within $\{(w,\xi):w\in\RR^d,\xi=\xi_i\}$ for $i\in\{1,2,...\}$.\\ 
	
	For an $a'\in\RR$, if $a'\geq a_j$,
	\begin{alignat}{6}
		&&&\{(w,\xi)\in\RR^{d+p}:\widetilde{\nabla}F_j(w,\xi)>a'\}&&=&\cup_{i=1}^{\infty}\{w\in\RR^d,\xi=\xi_i:\widetilde{\nabla}F_j(w,\xi)>a'\},\nonumber
	\end{alignat}
	otherwise
	\begin{alignat}{6}
		\{(w,\xi)\in\RR^{d+p}:\widetilde{\nabla}F_j(w,\xi)>a'\}
		&&=&\cup_{i=1}^{\infty}\{w\in\RR^d,\xi=\xi_i:\widetilde{\nabla}F_j(w,\xi)>a'\}\nonumber\\
		&&&\cup\{(w,\xi): w\in\RR^d,\xi\in\RR^p\setminus(\cup_{i=1}^{\infty}\xi_i)\},\nonumber
	\end{alignat}
	showing that $\widetilde{\nabla}F(w,\xi)$ is Borel measurable for $(w,\xi)\in\RR^{d+p}$.
\end{proof}

\begin{proof}{\it (Proposition \ref{omin})} For the goal of showing that functions are definable in the o-minimal structure of $\RR_{\text{exp},<}$, we will give a very short background on definable sets, which will serve our proofs. An {\bf atomic formula} is a relation symbol $\{>,=\}$ applied to {\bf terms} which are made up of finitely many applications of the functions $\{+,\cdot,exp\}$ to variables $\{w_1,w_2,...\}$ and constants taken from $\RR$.\footnote{O-minimality is shown for definable sets with parameters, in our case, taken from $\RR$, so we can extend the constants from $\{0,1\}$ to $\RR$.} {\bf Formulas} are finitely many applications of boolean operations $\{\lor,\land,\neg\}$ and the existential quantifier $\exists$ to atomic formulas. For a formula $\phi$, a definable set in $\RR^d$ is the subset $X\subseteq\RR^d$ such that $\phi$ is true. 
	Let $\Gamma_{g}:=\{(w,y)\in\RR^{d+q}:g(w)=y\}$ be the graph of a function $g:\RR^d\rightarrow \RR^q$. 
	A definable function is a function whose graph is definable, and a function $g:\RR^d\rightarrow \RR^q$ is definable if and only if each of its coordinate functions $g_i(w)$ for $i\in[q]$ are definable \citep[Exercise 1.10]{coste2000}. The composition of definable functions is definable, as is the addition and multiplication of definable functions \citep[Exercise 1.11]{coste2000}.\\  
	
	{\bf Affine Map (AM):}
	The graph of $\text{AM}(w,b):=\left<w,x\right>+b$ for variables $w\in\RR^d$, $b\in\RR$, and constants $x\in\RR^d$, 
	\begin{alignat}{6}
		\Gamma_{\text{AM}}=\{(w,b,y)\in\RR^{d+2}:\left<w,x\right>+b=y\}\nonumber
	\end{alignat}
	is definable, its gradient $\nabla \text{AM}(w,b)=[x^T,1]^T$ is measurable as a constant, with $\nabla \text{AM}(w,b)=\overline{\partial}\text{AM}(w,b)$ since it is continuous.\\
	
	{\bf ReLU:}
	The graph of $\text{ReLU}(w)$ for $w\in \RR$, 
	\begin{alignat}{6}
		\Gamma_{\text{ReLU}}=&\{(w,y)\in\RR^2:(w>0\land w=y)\lor(\neg(w>0)\land y=0)\}\nonumber
	\end{alignat} 
	is definable, its bp derivative $\text{ReLU}'(w)=\ind_{\{w\in\RR:w>0\}}$ is the indicator function of a measurable set, and for $w\in\RR$, $\text{ReLU}'(w)\in \overline{\partial} \text{ReLU}(w)$, which equals the subdifferential from convex analysis \citep[Proposition 2.2.7]{clarke1990}.\\
	
	{\bf Conv2d:} Each of its component functions is an affine mapping, i.e. a filter being applied to the input layer, hence each component function of Conv2d is definable with a measurable gradient.\\
	
	{\bf MaxPool2d:}
	Each of its component functions equals $\text{MP}(w):=\max\limits_{\substack{i\in [H]_{-1}\\j\in[W]_{-1}}}w_{ij}$ for a subset $w\in\RR^{H\times W}$ of the input layer. The graph of this function can be written as 
	\begin{alignat}{6}
		\Gamma_{\text{MP}}=\{(w,y)\in\RR^{H\times W+1}:\lor_{\substack{i\in [H]_{-1}\\j\in[W]_{-1}}}\big((w_{ij}=y)\land_{\substack{k\in [H]_{-1}\\l\in[W]_{-1}}}\neg(w_{kl}>w_{ij})\big)\}.\nonumber
	\end{alignat} 
	
	Looping through $(i,j)$, the bp gradient of $\text{MP}(w)$ is set to $1$ for the first pair of indices $(i,j)$ such that $w_{i,j}=\text{MP}(w)$, with the remaining entries of the gradient set to $0$. The bp gradient can be expressed recursively as
	\begin{alignat}{6}
		\nabla_{ij}\text{MP}(w)=\big(\prod\limits_{\substack{k\in[H]_{-1}\\l\in[W]_{-1}}}\ind_{\{w\in\RR^{H\times W}:w_{ij}\geq w_{kl}\}}\big)\big(1-\sum\limits_{\substack{k\in[i]_{-1}\\l=[W]_{-1}}}\nabla_{kl}\text{MP}(w)
		-\sum\limits_{l\in[j]_{-1}}\nabla_{il}\text{MP}(w)\big),\nonumber
	\end{alignat} 
	which is the product and subtraction of real-valued measurable functions. Let $E(w)$ be the set of pairs $(i,j)$ such that $w_{ij}=\text{MP}(w)$, $E(w):=\{(i,j)\in[H]_{-1}\times[W]_{-1}: w_{ij}=MP(w)\}$, and let $e_{ij}$ be a matrix of dimension $H\times W$ equal to 1 at entry $(i,j)$ and equal to $0$ otherwise, for each $(i,j)\in E(w)$. The Clarke subdifferential of $MP(w)$ equals $\overline{\partial}\text{MP}(w)=\co\{e_{ij}: (i,j)\in E(w)\}$ \citep[Proposition 2.3.12]{clarke1990}, hence $\nabla\text{MP}(w)\in\overline{\partial}\text{MP}(w)$.\\
	
	{\bf Crossentropyloss (CL):}
	For $C$ classes, Crossentropyloss takes the form of 
	\begin{alignat}{6}
		\text{CL}(w):=-\log\left(\frac{e^{w_t}}{\sum_{i=0}^{C-1}e^{w_i}}\right),\nonumber
	\end{alignat}
	where $t$ is the index of the target class. The graph of this function can be written as
	\begin{alignat}{6}
		\Gamma_{\text{CL}(w)}=\{(w,y)\in\RR^{C+1}:\sum_{i=0}^{C-1}e^{w_i}=e^ye^{w_t}\},\nonumber
	\end{alignat}
	hence $\text{CL}(w)$ is definable. The gradient is continuous, with components equal to 
	\begin{alignat}{6}
		\nabla_t CL(w)=\frac{e^{w_t}}{\sum_{i=0}^{C-1}e^{w_i}}-1\nonumber
	\end{alignat}
	and for $j\neq t$,
	\begin{alignat}{6}
		\nabla_j \text{CL}(w)=\frac{e^{w_j}}{\sum_{i=0}^{C-1}e^{w_i}},\nonumber
	\end{alignat}
	therefore measurable with $\nabla \text{CL}(w)=\overline{\partial}\text{CL}(w)$.
\end{proof} 

\section*{Appendix F. Details of Section \ref{exper}}

\subsubsection*{Details of the neural network architecture:}
Following Pytorch, let $\text{Conv2d}(i,o,k)$ denote a 2D convolutional layer with $i$ input and $o$ output channels, using $k\times k\times i$ sized filters, with a stride of 1 and 0 padding. Let $\text{MaxPool2d}(2,2)$ be a 2D max pool layer with a window size of $2\times 2$, stride of 2, and 0 padding, and let $\text{Linear}(i,o)$ be a fully connected layer with $i$ and $o$ being the number of inputs and outputs. The trained neural network then takes the following form:
\begin{alignat}{6}
	&\text{Input}\rightarrow \text{Conv2d}(1,6,5)\rightarrow \text{ReLu}\rightarrow\text{MaxPool2d}(2,2)
	\rightarrow\text{Conv2d}(6,16,5)\nonumber\\
	\rightarrow&\text{ReLu}
	\rightarrow\text{MaxPool2d}(2,2)
	\rightarrow\text{Conv2d}(16,120,4)\rightarrow\text{ReLu}\rightarrow\text{Linear}(120,84)\nonumber\\
	\rightarrow&\text{ReLu}\rightarrow\text{Linear}(84,10)	\rightarrow\text{CrossEntropyLoss}\rightarrow\text{Output}.\nonumber
\end{alignat}

\subsubsection*{Overview of the BNB implementation:}

We note that since $\{p_i\}\subset \ZZ_{>0}$ for this specific application, with $m=\lfloor (1-s)d\rfloor$, dynamic programming could have been used to compute $\Pi_{C\cap \overline{B}_{\beta}}(\cdot)$, but the BNB approach was implemented to be applicable for $\{p_i\}\subset \RR_{>0}$, $m\in\RR_{>0}$, and the real-valued objective coefficients of \eqref{knap2}, which follows Section \ref{constraint} where minimal assumptions were placed on $\{p_i\}$ and $m$.\\ 

For simplicity let $y_i:=(||w^i||^2_2-||\max(|w^i|-\beta,0)||^2_2)$. For each node of the search tree, where a 0-1 knapsack problem is considered with subsets of the decision variables already assigned 1 and 0, the lower bound of the problem uses the Greedy-Split algorithm described in \citep[Chapter 2.1]{kellerer2004} for the undetermined decision variables. Let $s$ be the index of the critical item \citep[Section 2.2]{kellerer2004}. The upper bound is computed following \citep[Equation 5.12]{kellerer2004} when there exists valid indices $s-1$ and $s+1$, or else as \citep[Equation 5.9]{kellerer2004} when $s=0$, where for both computations the floor operators are not employed as generally $\{y_i\}\notin\ZZ_{>0}$.\\

Before branching if the upper bound is greater than the global lower bound, significant algorithm speed up was observed by checking if $z_s$ dominates or is dominated by a $z_i$ already set to $0$ or $1$. We say that $z_j$ dominates $z_k$ if either $y_j\geq y_k$ and $p_j< p_k$ or $y_j>y_k$ and $p_j\leq p_k$. Before branching with $z_s=0$, we checked that $z_s$ does not dominate a $z_i=1$, and before branching with $z_s=1$, we checked that $z_s$ is not dominated by a $z_i=0$. If one of these cases occurred, the branch was abandoned as the resulting solution would not be optimal.

\subsubsection*{Estimation of $L_0$, $Q$, and $\Delta$:}
We took $q=250$ samples $\{\xi_i\}_{i=1}^q$ and $q$ pairs of sampled points $\{(w_j,v_j)\}_{j=1}^{q}\subset\overline{B}_{\kappa}$. The points $w_j$ are uniformly sampled in $\overline{B}_{\kappa-\iota}$ for $\iota=0.01/\kappa$, and each point $v_j$ is sampled uniformly near $w_j$, in $w_j+\overline{B}_{\iota}$. The estimate of $L_0(\xi_i)$ is 
\begin{alignat}{6}
	\widehat{L}_{0}(\xi_i)=\max\limits_{j\in[q]}\frac{|F(w_j)-F(v_j)|}{||w_j-v_j||_2}.\nonumber
\end{alignat}
The estimate of $L_0$ is $\widehat{L}_0=\mean(\widehat{L}_{0}(\xi_i))$ and the estimate of $Q$ equals $\widehat{Q}=\mean(\widehat{L}_{0}(\xi_i)^2)$.  For each run, an estimate of $\Delta$, $\overline{\Delta}_l$ for $l\in [3]$, was computed. Given the randomly generated $w_0^l$, $w^l_1=\Pi_{S\cap \overline{B}_{\beta}}(w^l_0)$. For each $\xi_i$ of the training set, $F(w_1^l+u_i,\xi_i)$ was computed for a sample $u_i\sim P_u$, and the average of these values, over $i$, were computed to estimate $f_{\alpha}(w_1^l)$, which was taken as $\overline{\Delta}_l$ given that $f_{\alpha}(w^*)\geq 0$. The estimate of $\Delta$ was then chosen as  $\overline{\Delta}=\max_{l\in[3]}\overline{\Delta}_l$. 

%\vskip 0.2in
\bibliography{sparse_training}

\end{document}